\documentclass[12pt,a4paper]{article}
\pdfoutput=1
\RequirePackage{amsmath}%
\usepackage{amsfonts}%
\usepackage{amssymb}%
\usepackage{graphicx}%
\usepackage{url}%
\usepackage{color}

\newtheorem{example}{Example}

\renewcommand{\Re}{\mathop{\mathrm{Re}}}
\renewcommand{\Im}{\mathop{\mathrm{Im}}}

\renewcommand{\i}{\mathrm{i}}

\newcommand{\bI}{{\bf I}}

\newcommand{\bM}{{\bf M}}
\newcommand{\bN}{{\bf N}}

\newcommand{\bn}{{\bf n}}

\begin{document}
\title{Simulating local fields in carbon nanotube reinforced composites for infinite strip with voids}
	\author{Mohamed Nasser$^{\rm a}$, El Mostafa Kalmoun$^{\rm b}$, Vladimir Mityushev$^{\rm c}$,\\ and Natalia Rylko$^{\rm c}$}
	
	\date{}
	\maketitle
	
	\vskip-0.8cm %
	\centerline{$^{\rm a}$Mathematics Program, Department of Mathematics, Statistics and Physics,} %
	\centerline{College of Arts and Sciences, Qatar University, Doha, Qatar}%
	%
	\centerline{$^{\rm b}$School of Science and Engineering, Al Akhawayn University in Ifrane,} %
	\centerline{PO Box 104, Ifrane 53000, Morocco}%
	%
	\centerline{$^{\rm c}$Faculty of Computer Science and Telecommunications, }
	\centerline{Cracow University of Technology, Krak\'{o}w, Poland}

\begin{abstract}
We consider the steady heat conduction problem within a thermal isotropic and homogeneous infinite strip composite reinforced by uniformly and randomly distributed non-overlapping carbon nanotubes (CNTs) and containing voids. We treat the CNTs as thin perfectly conducting elliptic inclusions and assume the voids to be of circular shape and act as barriers to heat flow.
We also impose isothermal conditions on the external boundaries by assuming the lower infinite wall to be a heater under a given temperature, and the upper wall to 
be 
a cooler that can be held at a lower fixed temperature. 
The equations for the temperature distribution are governed by the two-dimensional Laplace equation with mixed Dirichlet-Neumann boundary conditions. The resulting boundary value problem is solved using the boundary integral equation with the generalized Neumann kernel. We illustrate the performance of the proposed method through several numerical examples including the case of the presence a large number of CNTs and voids.
\end{abstract}

\begin{center}
\begin{quotation}
\end{quotation}
\end{center}


\section{Introduction}
Nanofibers  embedded in polymer matrices have attracted attention as one of the reinforcements for composite materials. Carbon nanotubes (CNTs) reinforced polymer nanocomposites are considered as conventional micro- and macro-composites \cite{Rahaman}. Their thermal, mechanical, and electric properties are determined by experimental and theoretical investigations \cite{Feng,Loos,Poveda}. CNTs are considered as perfectly conducting inclusions, which suggests imposing Dirichlet boundary conditions on the boundary of CNTs. On the other hand, the classical problems for materials with holes in porous media and materials with voids and insulting inclusions are modeled by the Neumann boundary condition \cite{PMA, GMN}. 

The present paper is devoted to the heat conduction within a 2D (two-dimensional) thermal isotropic and homogeneous nanocomposite, which takes the form of an infinite strip, when it is reinforced by non-overlapping and randomly distributed CNTs and contains defects and voids. In particular, we are interested in studying the effect of CNTs as well as of the presence of voids on the macroscopic conductive and mechanical properties of this composite. Owing to the superconductivity of CNTs and the extremely low conductivity of voids, we can assume that the conductivity of CNTs, of the polymer host and of voids to be governed by the inequalities $\lambda_c \gg \lambda \gg \lambda_0$. Such an assumption leads to a mixed boundary value problem where the latter inequality becomes $+\infty \gg \lambda \gg 0$. 
The host conductivity can be normalized to unity, i.e., $\lambda=1$.  

Theoretical investigation of mixed boundary value problems by integral equations can be found in~\cite{Nas-mix,Nas-Vou}. In the same time, implementation of  numerical methods for large number of inclusions and holes is still a challenging problem of applied and computational mathematics. We propose in this work a fast and effective algorithm for the numerical solution of the formulated mixed boundary value problem. The method is based on the boundary integral equation with the generalized Neumann kernel. 
The integral equation has been used in~\cite{Nas-Vou} to solve a similar mixed boundary value problem related to the capacity of generalized condensers. 
The proposed method can be even employed when the number of perfectly conducting inclusions and holes is very large.

As a result of simulations, we first study the 2D local fields for three types of media. In the first type, we consider the case of pure $m$ void cracks with $m=5, 30, 50$. The second type consists of pure $\ell$ CNT inclusions with $\ell=5$ and $200$. Finally, we treat the case of a large number of combined inclusions and holes by considering either $2000$ of one of the two or $1000$ of each. Afterward, we take up the systematic investigation of the effective conductivity of the considered composites. It is important in applications to predict the macroscopic properties of composites which depend on the concentration of perfectly conducting CNTs as well as on the concentration of holes and voids. It is worth noting that the notation of concentration are different for slit shapes of CNTs and circular shapes of holes.  The performed simulations of local fields and computation of their averaged conductivities for various concentrations allows to establish the dependence of the macroscopic conductivity on the main geometrical parameters.

\section{Problem formulation}

\begin{figure}[ht] %
\centerline{
\includegraphics[page=1, trim =0cm 0cm 0cm 0cm, clip, width=0.7\textwidth]{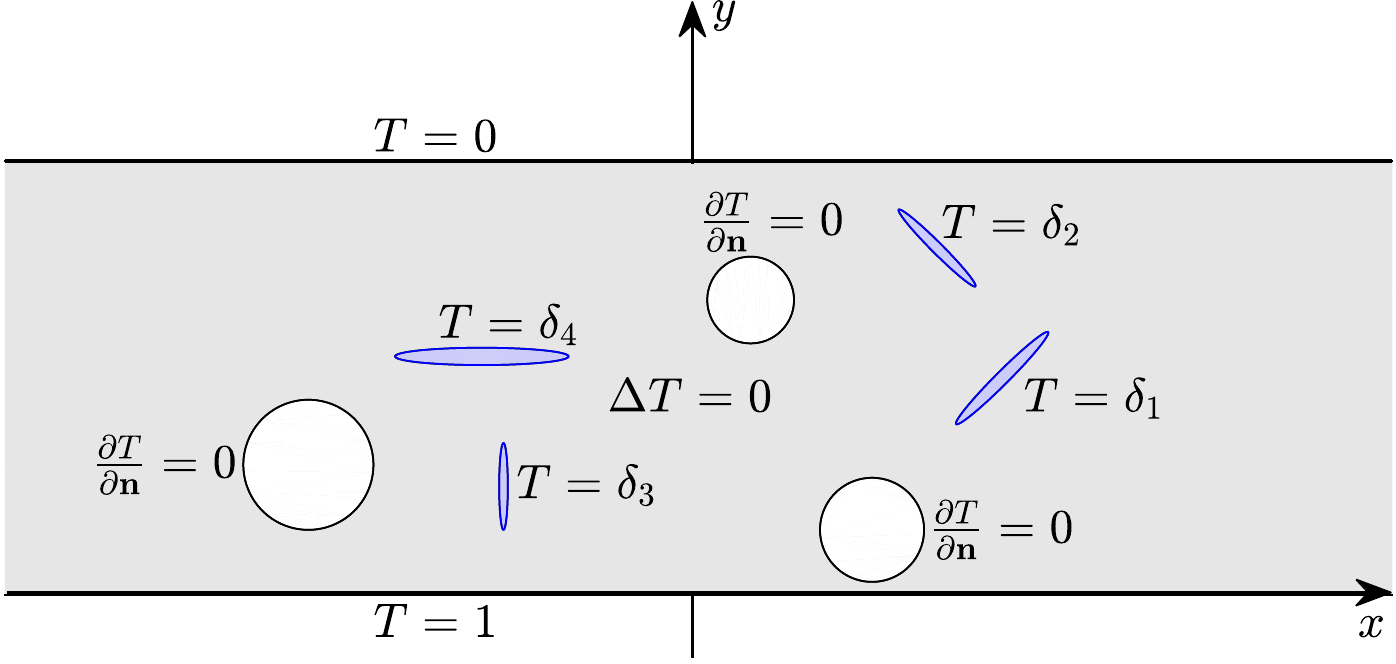}
}
\caption{Geometry of the problem (for $\ell=4$ and $m=7$).}
\label{fig:dom-Om}
\end{figure}

Let us consider a channel medium embedding $m$ inhomogeneities in the form of $\ell$ nanofillers and $p=m-\ell$ holes (voids). As many nanofillers (e.g, carbon nanontubes~\cite{malekie2015study}) have cross sections of elliptical shapes, we model them as ellipses $C_{1},\ldots, C_{\ell}$. Furthermore, we represent the non-conducting holes by inner circles 
 $C_{\ell+1},\ldots, C_m$. The top and bottom infinite walls of the channel are denoted respectively by $C_0'$ and $C_0''$, which yields a multiply connected domain $\Omega $ of connectivity $m+1$ with a boundary set $C =\bigcup_{k=0}^{m}C_k$ where $C_0=C_0'\cup C_0''$. An example of this domain for the case of $\ell=4$ and $m=7$ is illustrated in Figure~\ref{fig:dom-Om}.

The medium matrix without inhomogeneities is supposed to be homogeneous and isotropic with a constant thermal conductivity $\lambda =1$.  We also assume that conduction is the only dominating mechanism of heat transfer in the medium.
Except being non-overlapping, no other restriction is imposed on the inhomogeneities as they can be placed at random orientation and position. 

The nanofillers are treated as heat superconductors with an almost uniform temperature distribution within each one. 
Therefore the temperature $T$ is assumed to be fixed to an indeterminate constant value $\delta_k$ along each ellipse $C_k$ for $k=1\ldots,\ell$. This assumption is consistent with the numerical results reported in~\cite{zhang2004simplified} for CNT reinforced polymer composites.
Furthermore, by the law of energy conservation in steady-state heat conduction, there should be no net thermal flow through each nanofiller. This constraint is written by means of the net heat flux boundary condition~\eqref{eq:bvp-nq}.

On the other hand, the curves $C_{\ell+1},\ldots, C_m$ are assumed to be perfect insulators and therefore they act as barriers to heat flow. Henceforth,  the Neumann boundary condition~\eqref{eq:bvp-q} is imposed along the holes contours.  
  Finally, isothermal conditions are imposed on the external boundaries  by assuming that the lower infinite wall is a heater of temperature $T_1$, and the upper wall acts as a heat sink, which can be held at a fixed temperature $T_0<T_1$. Thee two values $T_0$ and $T_1$ of the temperature on the external boundaries are normalized to $0$ and $1$, respectively.
 
Under steady-state conditions, Fourier's law of heat conduction and the above specified heat boundaries conditions yield the temperature distribution $T$ governed by the following mixed Dirichlet-Neumann boundary value problem: 
\begin{subequations}\label{eq:mix-bd-T}
	\begin{align}
	\label{eq:bvp-1-T}
	\Delta T &= 0 \quad \mbox{in }\Omega, \\
	\label{eq:bvp-0-Ta}
	T&= 0 \quad \mbox{on }C'_{0}, \\
	\label{eq:bvp-0-Tb}
	T&= 1 \quad \mbox{on }C''_{0}, \\
	\label{eq:bvp-3-T}
	T &= \delta_k \quad \mbox{on }C_{k}, \quad k=1,2,\ldots,\ell, \\
	\label{eq:bvp-nq}
	\int_{C_{k}}\frac{\partial T}{\partial\bn}ds &= 0 \quad k=1,2,\ldots,\ell, \\
	\label{eq:bvp-q}
	\frac{\partial T}{\partial\bn}&= 0 \quad \mbox{on }C_{k}, \quad k=\ell+1,\ell+2,\ldots,m,
	\end{align}
\end{subequations}
where $\partial T/\partial\bn$ denotes the normal derivative of $T$, and $\delta_1,\ldots,\delta_m$ are undetermined real constants that need to be found alongside the distribution temperature $T$.

\section{The integral equation method}

\begin{figure}[ht] %
\centerline{
\includegraphics[page=1, trim =0cm 0cm 0cm 0cm, clip, width=0.4\textwidth]{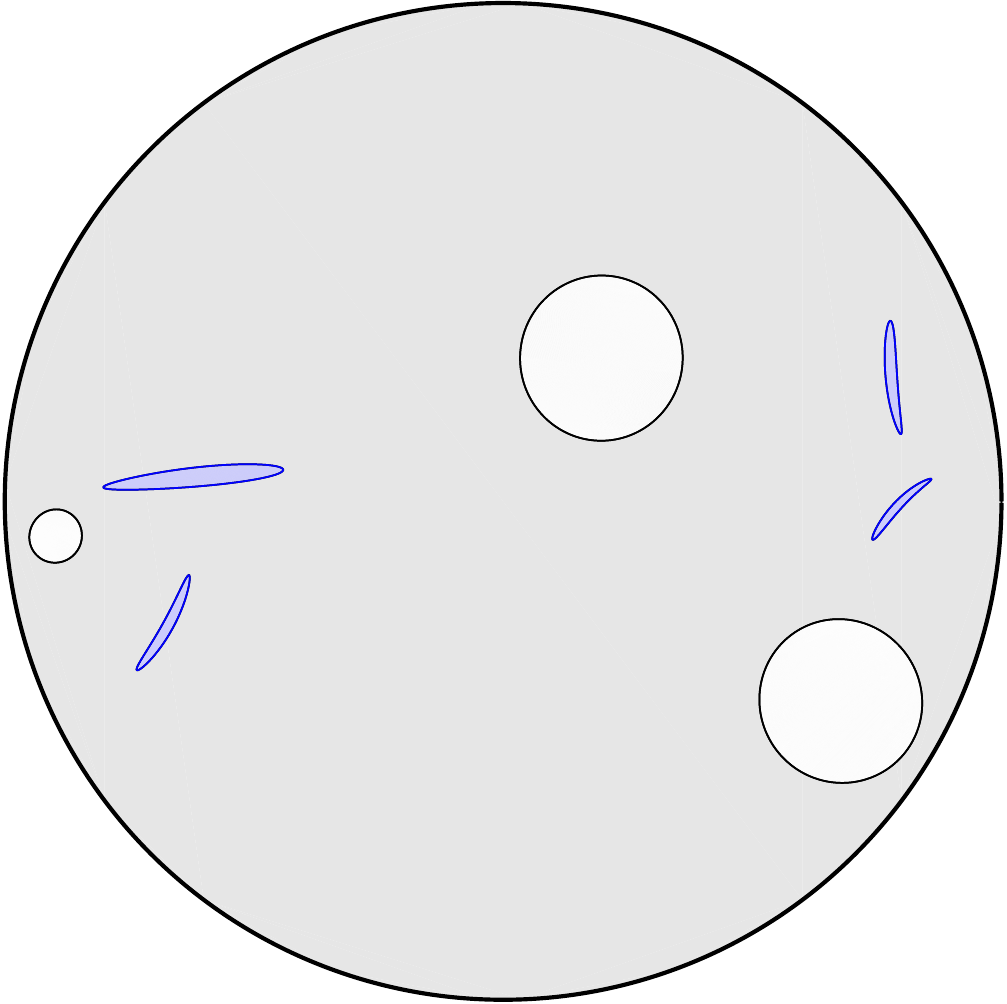}
}
\caption{The computational domain $G$ corresponding to the physical domain in Figure~\ref{fig:dom-Om}.}
\label{fig:dom-G}
\end{figure}

The boundary integral equation with the generalized Neumann kernel is not directly applicable to the above boundary value problem~(\ref{eq:mix-bd-T}) because of the external boundary component. However, the boundary value problem~(\ref{eq:mix-bd-T}) is invariant under conformal mapping. The mapping function
\[
z = \Phi(\zeta) = \frac{1}{\pi}\log\frac{1+\zeta}{1-\zeta}+\frac{\i}{2}
\]
conformally maps the unit disk $|\zeta|<1$ onto the infinite strip $0<\Im z<1$. Thus, the inverse mapping
\[
\zeta = \Phi^{-1}(z) = \tanh\left(\frac{\pi z}{2}-\frac{\pi\i}{4}\right)
\]
conformally maps the infinite strip $0<\Im z<1$ onto the unit disk $|\zeta|<1$, the real axis onto the lower half of the unit circle, the line $\Im z=1$ onto  the upper half of the unit circle, and satisfies $\Phi^{-1}(\pm\infty+0\i)=\pm1$. Consequently, the function $\Phi^{-1}$ maps the multiply connected domain $\Omega$ in the $z$-plane (the physical domain) onto a multiply connected domain $G$  in the $\zeta$-plane interior of the unit circle and exterior of $m$ smooth Jordan curves (the computational domain). In Figure~\ref{fig:dom-G}, we display the result of the conformal mapping of the example shown in Figure~\ref{fig:dom-Om}.

It follows that the harmonic function $T$ can be written as 
\[
T(z)=U(\Phi^{-1}(z))
\]
 in which the function $U$ is the solution of the following boundary value problem in the $\zeta$-plane:
\begin{subequations}\label{eq:mix-bd-U}
	\begin{align}
	\label{eq:bvp-1-U}
	\Delta U &= 0 \quad \mbox{in }G, \\
	\label{eq:bvp-0-Ua}
	U&= 0 \quad \mbox{on }\Gamma'_{0}, \\
	\label{eq:bvp-0-Ub}
	U&= 1 \quad \mbox{on }\Gamma''_{0}, \\
	\label{eq:bvp-3-U}
	U &= \delta_k \quad \mbox{on }\Gamma_{k}, \quad k=1,2,\ldots,\ell, \\
	\label{eq:bvp-4-U}
	\int_{\Gamma_{k}}\frac{\partial U}{\partial\bn}ds &= 0 \quad k=1,2,\ldots,\ell, \\
	\label{eq:bvp-2-U}
	\frac{\partial U}{\partial\bn}&= 0 \quad \mbox{on }\Gamma_{k}, \quad k=\ell+1,\ell+2,\ldots,m,
	\end{align}
\end{subequations}
where $\Gamma'_{0}=\Phi^{-1}(C'_{0})$, $\Gamma''_{0}=\Phi^{-1}(C''_{0})$, and $\Gamma_{k}=\Phi^{-1}(C_{k})$ for $k=1,2,\ldots,m$.
Note that the restriction of the function $U(\zeta)$ on the external boundary is discontinuous at $\zeta=\pm1$. However, the function $U$ can be cast into the form
\[
U(\zeta)=u_0(\zeta)+ u(\zeta)
\]
where $u(\zeta)$ is a harmonic function in $G$, and
\[
u_0(\zeta)=\frac{1}{\pi}\Im\log\frac{1-\zeta}{1+\zeta}+\frac{1}{2}.
\]
The function $u_0(\zeta)$ is harmonic in $G$ with $u_0(\zeta)=0$ on the upper half of the unit circle and $u_0(\zeta)=1$ on the lower part. The function $u(\zeta)$ is the solution of the boundary value problem
\begin{subequations}\label{eq:mix-bd-ut}
	\begin{align}
	\label{eq:bvp-1-u}
	\Delta u(\zeta) &= 0 \quad \mbox{if }\zeta\in G, \\
	\label{eq:bvp-0-u}
	u(\zeta)&= 0 \quad \mbox{if }\zeta\in\Gamma_{0}, \\
	\label{eq:bvp-3-u}
	u(\zeta) &= \delta_k-\frac{1}{\pi}\Im\log\frac{1-\zeta}{1+\zeta}-\frac{1}{2} \quad \mbox{if }\zeta\in\Gamma_{k}, \quad k=1,2,\ldots,\ell, \\
	\label{eq:bvp-4-u}
	\int_{\Gamma_{k}}\frac{\partial u}{\partial\bn}ds &= 0 \quad k=1,2,\ldots,\ell, \\
	\label{eq:bvp-2-u}
	\left.\frac{\partial u}{\partial\bn}\right|_{\zeta}&= -\left.\frac{\partial u_0}{\partial\bn}\right|_{\zeta} \quad \mbox{if }\zeta\in\Gamma_{k}, \quad k=\ell+1,\ell+2,\ldots,m,
	\end{align}
\end{subequations}
where $\Gamma_0$ is the unit circle.

For the orientation of the boundary components of $G$, we assume that $\Gamma_0$ is oriented counterclockwise and the other curves $\Gamma_1,\ldots,\Gamma_m$ are oriented clockwise. We assume that each boundary component $\Gamma_k$, $k=0,1,\ldots,m$, is parametrized by a $2\pi$-periodic function $\eta_k(t)$, $t\in J_k:=[0,2\pi]$ such that $\eta'_k(t)\ne0$. Let $J$ be the disjoint union of the $m+1$ intervals $J_0,\ldots,J_m$, the whole boundary $\Gamma$ is parametrized by the complex function $\eta$ defined on $J$ by~\cite{Weg-Nas,Nas-ETNA}
\[
\eta(t)= \left\{ \begin{array}{l@{\hspace{0.5cm}}l}
	\eta_0(t),&t\in J_0,\\
	\eta_1(t),&t\in J_1,\\
	\hspace{0.3cm}\vdots\\
	\eta_m(t),&t\in J_m.
	\end{array}
	\right.
\]
Note that the unit circle $\Gamma_0$ is parametrized by $\eta_0(t)=e^{\i t}$, $t\in J_0=[0,2\pi]$.

Let $\bn(\zeta)$ be the unit outward normal vector at $\zeta\in\Gamma$ and let $\nu(\zeta)$ be the angle between the normal vector $\bn(\zeta)$ and the positive real axis. Then, for $\zeta=\eta(t)\in\Gamma$, 
\begin{equation}\label{eq:n-nu}
\bn(\zeta) = e^{\i\nu(\zeta)}=-\i\frac{\eta'(t)}{|\eta'(t)|}.
\end{equation}
Thus
\begin{equation}\label{eq:nd-u0}
\frac{\partial u_0}{\partial \bn}=\nabla u_0\cdot\bn
=\cos\nu\frac{\partial u_0}{\partial x}+\sin\nu\frac{\partial u_0}{\partial y}
=\Re\left[e^{\i\nu}
\left(\frac{\partial u_0}{\partial x}-\i\frac{\partial u_0}{\partial y}\right)\right].
\end{equation}

The harmonic function $u_0(\zeta)$ is the real part of a single-valued analytic function $f_0(\zeta)$, i.e., $u_0(\zeta)=\Re[f_0(\zeta)]$, where
\begin{equation}\label{eq:f-0}
f_0(\zeta) = \frac{1}{\pi\i}\log\frac{1-\zeta}{1+\zeta}+\frac{1}{2},
\end{equation}
and the branch of the logarithm function is chosen such that $\log 1=0$. Then by the Cauchy-Riemann equations, we have
\[
f'_0(\zeta)=\frac{\partial u_0(\zeta)}{\partial x}-\i\frac{\partial u_0(\zeta)}{\partial y},
\]
which, in view of~\eqref{eq:n-nu} and~\eqref{eq:nd-u0}, implies that
\begin{equation}\label{eq:u-normal-0}
|\eta'(t)|\,\left.\frac{\partial u_0}{\partial \bn}\right|_{\eta(t)}
=\Re\left[-\i\eta'(t)\,f'_0(\eta(t))\right], \quad \eta(t)\in\Gamma_k, \quad k=\ell+1,\ldots,m.
\end{equation}
Since 
\[
f'_0(\zeta) = \frac{\i}{\pi}\left(\frac{1}{1-\zeta}+\frac{1}{1+\zeta}\right),
\]
it follows that for $\eta(t)\in\Gamma_k$ and $k=\ell+1,\ldots,m$,
\begin{equation}\label{eq:u0-normal}
|\eta'(t)|\,\left.\frac{\partial u_0}{\partial \bn}\right|_{\zeta=\eta(t)}
=\frac{1}{\pi}\Re\left[\frac{\eta'(t)}{1-\eta(t)}+\frac{\eta'(t)}{1+\eta(t)}\right].
\end{equation}

The harmonic function $u$ can be assumed to be a real part of an analytic function $f(\zeta)$, $\zeta\in G$.
The boundary conditions~(\ref{eq:bvp-0-u}) and~(\ref{eq:bvp-3-u}) give the real parts of the function $f(\zeta)$  on $\Gamma_k$ for $k=0,1,\ldots,\ell$. Specifically, we have
\begin{equation}\label{eq:delta-k-1}
\Re[f(\eta(t))]=0, \quad \eta(t)\in\Gamma_0
\end{equation}
and
\begin{equation}\label{eq:delta-k-2}
\Re[f(\eta(t))]=\delta_k-\frac{1}{\pi}\Im\log\frac{1-\eta(t)}{1+\eta(t)}-\frac{1}{2} \quad \mbox{if }\eta(t)\in\Gamma_{k},  \quad k=1,\ldots,\ell.
\end{equation}
For the remaining boundary components $\Gamma_k$ for $k=\ell+1,\ldots,m$, we use the condition~(\ref{eq:bvp-2-u}) to determine the boundary condition on $f(\eta)$. 
By the Cauchy-Riemann equations, we can show using the same arguments as in~(\ref{eq:u-normal-0}) that
\begin{equation}\label{eq:u-normal}
|\eta'(t)|\,\left.\frac{\partial u}{\partial \bn}\right|_{\eta(t)}
=\Re\left[-\i\eta'(t)\,f'(\eta(t))\right]
\end{equation}
Thus, for $\eta(t)\in\Gamma_k$ for $k=\ell+1,\ldots,m$, it follows from~(\ref{eq:bvp-2-u}), (\ref{eq:u0-normal}), and~(\ref{eq:u-normal}) that
\[
\Re\left[-\i\eta'(t)\,f'(\eta(t))\right]
=-\frac{1}{\pi}\Re\left[\frac{\eta'(t)}{1-\eta(t)}+\frac{\eta'(t)}{1+\eta(t)}\right].
\]
Integrating with respect to the parameter $t$ for $t\in J_k$, $k=\ell+1,\ldots,m$, we obtain
\begin{equation}\label{eq:delta-k-3}
\Re\left[-\i f(\eta(t))\right]
=\frac{1}{\pi}\log\left|\frac{1-\eta(t)}{1+\eta(t)}\right|+\delta_k,
\end{equation}
where the integration constants $\delta_k$ are undetermined. 
The constants $\delta_k$, $k=1,\ldots,m$ in~(\ref{eq:delta-k-2}) and~(\ref{eq:delta-k-3}) are determined so that $f(z)$ is a single-valued analytic function. 

Since we are interested in the function $u$, the real part of $f$, we may assume that $c=f(\alpha)$ is real for some given point $\alpha$ in $G$.
Define an analytic function $g(\zeta)$ in the domain $G$ through
\begin{equation}\label{eq:f-g}
f(\zeta)=(\zeta-\alpha)g(\zeta)+c.
\end{equation}
Define also
\begin{equation}\label{eq:A}
A(t)=e^{-\i\theta(t)}(\eta(t)-\alpha),
\end{equation}
where $\theta(t)$ is the piecewise constant function given by
\begin{equation}\label{eq:thet}
\theta(t)=\left\{ \begin{array}{l@{\hspace{0.5cm}}l}
0,&t\in J_0,\\
\hspace{0.3cm}\vdots\\
0,&t\in J_\ell,\\
\pi/2,&t\in J_{\ell+1},\\
\hspace{0.3cm}\vdots\\
\pi/2,&t\in J_{m}.
\end{array}
\right.
\end{equation}
Thus
\[
e^{-\i\theta(t)}f(\eta(t))=A(t)g(\eta(t))+e^{-\i\theta(t)}c,
\]
which implies that
\[
\Re[A(t)g(\eta(t))]=\Re[e^{-\i\theta(t)}f(\eta(t))]-c\cos\theta(t).
\]

On the basis of the conditions~(\ref{eq:delta-k-1}), (\ref{eq:delta-k-2}), and~(\ref{eq:delta-k-3}), the function $g(z)$ satisfies the Riemann-Hilbert problem
\begin{equation}\label{eq:rhp}
\Re[A(t)g(\eta(t))]=\gamma(t)+h(t),
\end{equation}
where
\begin{equation}\label{eq:gam}
h(t)=\left\{ \begin{array}{l@{\hspace{0.5cm}}l}
-c,&t\in J_0,\\
\delta_1-\frac{1}{2}-c,&t\in J_1,\\
\hspace{0.3cm}\vdots\\
\delta_\ell-\frac{1}{2}-c,&t\in J_\ell,\\
\delta_{\ell+1},&t\in J_{\ell+1},\\
\hspace{0.3cm}\vdots\\
\delta_m,&t\in J_{l+p}.
\end{array}
\right.,\quad
\gamma(t)=\left\{ \begin{array}{l@{\hspace{0.5cm}}l}
0,&t\in J_0,\\
-\frac{1}{\pi}\Im\log\frac{1-\eta(t)}{1+\eta(t)},&t\in J_1,\\
\hspace{0.3cm}\vdots\\
-\frac{1}{\pi}\Im\log\frac{1-\eta(t)}{1+\eta(t)},&t\in J_l,\\
\frac{1}{\pi}\log\left|\frac{1-\eta(t)}{1+\eta(t)}\right|,&t\in J_{\ell+1},\\
\hspace{0.3cm}\vdots\\
\frac{1}{\pi}\log\left|\frac{1-\eta(t)}{1+\eta(t)}\right|,&t\in J_{m}.
\end{array}
\right.
\end{equation}
It is clear that the function $\gamma$ is known and the piecewise constant function $h$ is unknown and should be determined. Let $\mu(t)=\Im[A(t)g(\eta(t))]$, i.e., the  boundary values of an analytic function $g$ are given by
\begin{equation}\label{eq:get}
g(\eta(t))=\frac{\gamma(t)+h(t)+\i\mu(t)}{A(t)}, \quad t \in J.
\end{equation}
Thus, in order to find the boundary values of the analytic function $g$, we need to determine the two unknown functions $\mu$ and $h$. These two functions can be computed using the boundary integral equation with the generalized Neumann kernel~\cite{Weg-Nas,Nas-ETNA,Nas-JMAA1}.

Let $H$ be the space of all real H\"older continuous functions on $\Gamma$, let $\bI$ be the identity operator, and let the integral operators $\bN$ and $\bM$ are defined on $H$ by
\begin{eqnarray*}
\bN\mu(s) &=& \int_J \frac{1}{\pi}\Im\left(
\frac{A(s)}{A(t)}\frac{\eta'(t)}{\eta(t)-\eta(s)}\right) \mu(t) dt, \quad s\in J,\\
\bM\mu(s) &=& \int_J \frac{1}{\pi}\Re\left(
\frac{A(s)}{A(t)}\frac{\eta'(t)}{\eta(t)-\eta(s)}\right) \mu(t) dt, \quad s\in J.
\end{eqnarray*}
The kernel of the operator $\bN$ is known as the generalized Neumann kernel. For more details, see~\cite{Weg-Nas,Nas-ETNA,Nas-JMAA1}. 
On account of~\cite{Nas-JMAA1}, we have $\mu$ is the unique solution of the integral equation
\begin{equation}\label{eq:ie}
(\bI-\bN)\mu=-\bM\gamma.
\end{equation}
Additionally, the piecewise constant function $h$ is given by
\begin{equation}\label{eq:h}
h=[\bM\mu-(\bI-\bN)\gamma]/2.
\end{equation}

We compute approximations to the functions $\mu$ in~\eqref{eq:ie} and $h$ in~\eqref{eq:h} by the MATLAB function \verb|fbie| from~\cite{Nas-ETNA}. This function employs a discretization of the integral equation~(\ref{eq:ie}) by the Nystr\"om method using the trapezoidal rule~\cite{Atk97} to obtain an algebraic linear system of size $(m+1)n\times(m+1)n$ where $n$ is the number of discretization points in each boundary component. The resulting system is solved  
by applying the generalized minimal residual method through  
the MATLAB function $\mathtt{gmres}$. The matrix-vector multiplication in $\mathtt{gmres}$ is computed using the MATLAB function $\mathtt{zfmm2dpart}$ from the $\mathtt{FMMLIB2D}$ MATLAB toolbox~\cite{Gre-Gim12}. The values of the other parameters in the function \verb|fbie| are chosen as in~\cite{NasMo}. For more details, we refer the reader to~\cite{Nas-ETNA}.

\section{Computing the temperature distribution and the heat flux}
\label{sc:q}

By computing $\mu$ and $h$, we obtain the boundary values of the function $g$ through~\eqref{eq:get}. 
The values of the function $g(\zeta)$ for $\zeta\in G$ can be computed by the Cauchy integral formula. For the numerical computation of $g(\zeta)$ for $\zeta\in G$, we use the MATLAB function \verb|fcau| from~\cite{Nas-ETNA}.
Then, the values of $f(\zeta)$ can be computed by~\eqref{eq:f-g} and hence the values of the solution of the boundary value problem~\eqref{eq:mix-bd-U} is given for $\zeta\in G$ by
\[
U(\zeta) = \Re\left[f(\zeta)+f_0(\zeta)\right].
\]
We deduce the values of the temperature distribution $T(z)$  for any $z\in \Omega$ by
\[
T(z) = \Re\left[f(\Phi^{-1}(z))+f_0(\Phi^{-1}(z))\right].
\]
Moreover, by computing the piecewise constant function $h$, we can compute as well the values of the undetermined real constants $c,\delta_1,\ldots,\delta_m$ from~\eqref{eq:rhp}.

The function $T(z)$ is the real part of the function
\[
F(z)=f(\Phi^{-1}(z))+f_0(\Phi^{-1}(z)), \quad z\in \Omega.
\]
According to the Cauchy-Riemann equations, it follows that the derivative of the complex potential $F(z)$ on $\Omega$ is given by
\[
F'(z)=\frac{\partial T}{\partial x}-\i\frac{\partial T}{\partial y}. 
\]
One the other hand,
\begin{equation}\label{eq:Fp}
F'(z)=\frac{f'(\Phi^{-1}(z))}{\Phi'(\Phi^{-1}(z))}+
\frac{f'_0(\Phi^{-1}(z))}{\Phi'(\Phi^{-1}(z))}, \quad z\in \Omega,
\end{equation}
where the denominator does not vanish in the domain $\Omega$ since $\Phi$ is a conformal mapping.
Therefore the heat flux can be expressed for $z\in D$ in terms of $F'(z)$ by the formula
\begin{equation}\label{eq:q}
q(z)=-\left.\left(\frac{\partial T}{\partial x},\frac{\partial T}{\partial y}\right)\right|_z
\equiv -\overline{F'(z)}.
\end{equation}
Hence
\begin{equation}\label{eq:dT-dy}
\frac{\partial T}{\partial y} = -\Im F'(z).
\end{equation}

The derivatives $f'_0(\Phi^{-1}(z)$ and $\Phi'(\Phi^{-1}(z))$ in~\eqref{eq:Fp} can be computed analytically. So, the values of the heat flux $q$ can be estimated on the domain $\Omega$ by first approximating the derivatives of the boundary values of the analytic function $f$ on each boundary components. This can be done by approximating the function $f(\eta(t))$ using trigonometric interpolating polynomials then differentiating. The values of $f'(\Phi^{-1}(z))$, in the right-hand side of~\eqref{eq:Fp}, can be then computed for $z\in \Omega$ using the Cauchy integral formula.

\section{Computing the effective thermal conductivity}
\label{sc:cond}

The medium matrix without inhomogeneities is assumed to be homogeneous and isotropic. 
We will assume that the CNTs and the circular voids are in the part of the domain between $x=-1$ and $x=1$.
Thus, the effective conductivity of a layer in the $y$-direction $\lambda_y$ is calculated by the formula (3.2.33) from the book \cite[p.~53]{GMN}, which in our case on account of~\eqref{eq:dT-dy} becomes
\begin{equation}
\lambda_y = -\frac 12 \int_{-1}^1 \frac{\partial T}{\partial y}(x,0) \;\textrm{d}x
=
\label{eq:lam_y}
\frac{1}{2}\Im\left[\int_{-1}^1 F'(x) \;\textrm{d}x\right].
\end{equation}

Since 
\[
\xi_0(t)=\Phi(\eta_0(t))=\Phi(e^{\i t})
=\frac{1}{\pi}\log\frac{1+e^{\i t}}{1-e^{\i t}}+\frac{\i}{2}, 
\quad 0\le t\le 2\pi, 
\]
where for $0<t<\pi$, $\xi_0(t)$ is on the line $y=1$ and for $\pi<t<2\pi$, $\xi_0(t)$ is on the real line. Thus, for $\pi<t<2\pi$, we have
\[
\xi_0(t)=\frac{1}{\pi}\log\left|\cot\frac{t}{2}\right|. 
\]
Since $-1=\xi_0(t_1)$ and $1=\xi_0(t_2)$ where
\begin{equation}\label{eq:t1t2}
\pi<t_1=2\pi-2\tan^{-1}\left(e^\pi\right)<t_2=2\pi-2\tan^{-1}\left(e^{-\pi}\right)<2\pi.
\end{equation}
Consequently, \eqref{eq:lam_y} can be written as
\begin{equation}\label{eq:lam_y2}
\lambda_y = \frac{1}{2}\Im\left[\int_{t_1}^{t_2} F'(\xi_0(t)) \xi'_0(t)\;\textrm{d}t\right].
\end{equation}
In combining~\eqref{eq:Fp} with the fact that $\xi'_0(t)=\i e^{\i t}\Phi'(e^{\i t})$, we can see that 
\[
F'(\xi_0(t))=\frac{f'(\Phi^{-1}(\xi_0(t)))}{\Phi'(\Phi^{-1}(\xi_0(t)))}+
\frac{f'_0(\Phi^{-1}(\xi_0(t)))}{\Phi'(\Phi^{-1}(\xi_0(t)))}
=\frac{f'(e^{\i t})}{\Phi'(e^{\i t})}+
\frac{f'_0(e^{\i t})}{\Phi'(e^{\i t})}. 
\]
Hence, 
\begin{equation}\label{eq:lam_y3}
\lambda_y = \frac{1}{2}\Im\left[\int_{t_1}^{t_2} \left[\i e^{\i t}\left(f'(e^{\i t})+f'_0(e^{\i t})\right)\right]\textrm{d}t\right],
\end{equation}
which implies that
\begin{equation}\label{eq:lam_y4}
\lambda_y = \frac{1}{2}\Im\left[f(e^{\i t_2})-f(e^{\i t_1})\right]+\frac{1}{2}\Im\left[f_0(e^{\i t_2})-f_0(e^{\i t_1})\right].
\end{equation}
The second term in the right-hand side of~\eqref{eq:lam_y3} does not depend on the CNTs or the voids. In view of~\eqref{eq:f-0} and~\eqref{eq:t1t2}, we have
\[
\frac{1}{2}\Im\left[f_0(e^{\i t_2})-f_0(e^{\i t_1})\right]=1,
\]
and hence
\begin{equation}\label{eq:lam_y5}
\lambda_y = 1+\frac{1}{2}\Im\left[f(e^{\i t_2})-f(e^{\i t_1})\right].
\end{equation}
Since $e^{\i t_1}$ and $e^{\i t_2}$ are on the unit circle $\Gamma_0$, the external boundary of $G$, and taking into account~\eqref{eq:f-g}, \eqref{eq:A}, \eqref{eq:thet}, and~\eqref{eq:get}, Equation~\eqref{eq:lam_y5} can be written as
\begin{equation}\label{eq:lam_y6}
\lambda_y = 1+\frac{1}{2}\left[\mu(t_2)-\mu(t_1)\right].
\end{equation}

By solving the integral equation~\eqref{eq:ie}, we obtain approximate values of $\mu$ at the discretization points. These values are employed to interpolate the approximate solution $\mu$ on $J_0$ by a trigonometric interpolation polynomial, which is then used to approximate the values of $\mu(t_1)$ and $\mu(t_2)$.					

\section{Numerical results}
\label{sc:num}

The above proposed method with $n=2^{11}$ is applied to compute the temperature field $T$ and the heat flux $q$ for several examples.
We will choose the CNTs and the circular voids within the part of the domain between $x=-1$ and $x=1$. 
To compute the values of the temperature distribution $T$ and the heat flux $q$, we discretize part of the domain $\Omega$, namely for $-1.5\le x\le 1.5$ and $0.0001\le y\le 0.9999$. Afterwards, we compute the values of the temperature distribution $T$ and the heat flux $q$ at these points as described in Section~\ref{sc:q}.

\subsection{The domain $\Omega$ with only circular voids}

In this subsection, we consider the domain $\Omega$ with $m$ non-overlapping circular holes and without any CNT (i.e., $\ell=0$). We also assume that all circular holes have the same radius $r$ with the parametrization
\[
\eta_j(t) = z_j+r e^{-\i t}, \quad 0\le t\le 2\pi, \quad j=1,2,\ldots,m,
\]
where $z_1,z_2,\ldots,z_m$ are the centers of the circular holes. As these circular holes are chosen in the part of the domain $\Omega$ between $x=-1$ and $x=1$, we  define the concentration $c(m,r)$ of these voids to be the area of these circular holes over the area of the rectangle $\{(x,y)\,:\,-1\le x\le 1,\;0\le y\le 1\}$, i.e.,  
\begin{equation}\label{eq:cir-conc}
c(m,r) = \frac{m r^2\pi}{2}.
\end{equation}

The Clausius-Mossotti approximation (CMA) also known as Maxwel's formula can be applied for dilute composites when the concentration \eqref{eq:cir-conc} is sufficiently small. Below, we write this formula for a macroscopically isotropic media with insulators of  identical circular holes within the precisely established precision in \cite{MR}
\begin{equation}\label{eq:CMA}
\lambda_e = \frac{1-c}{1+c} + O(c^3).
\end{equation}

\begin{figure}[ht] %
\centerline{
\includegraphics[page=1, trim =0.65cm 2.25cm 0.65cm 2.75cm, clip, width=0.52\textwidth]{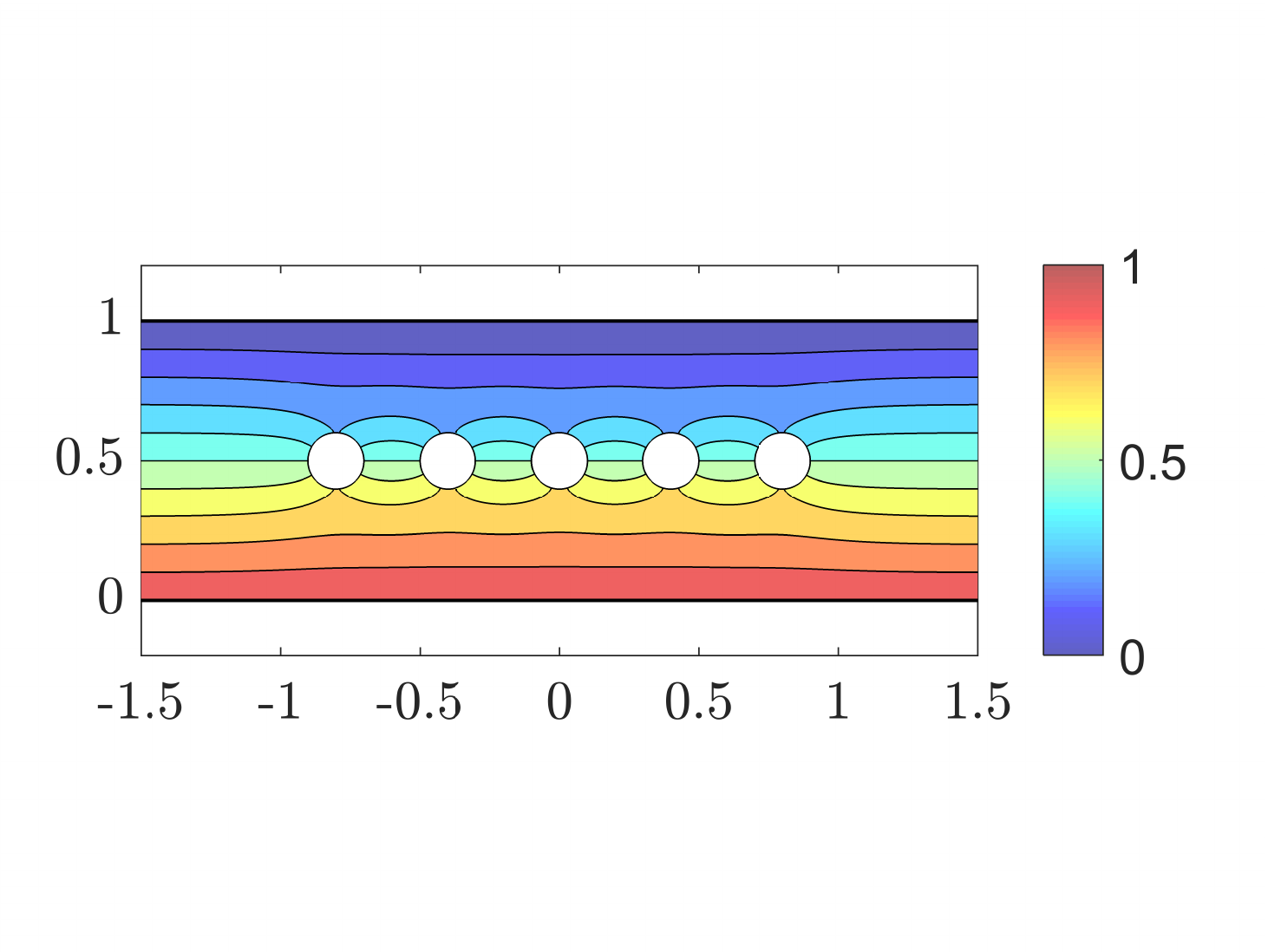}
\hfill
\includegraphics[page=1, trim =0.65cm 2.25cm 0.65cm 2.75cm, clip, width=0.52\textwidth]{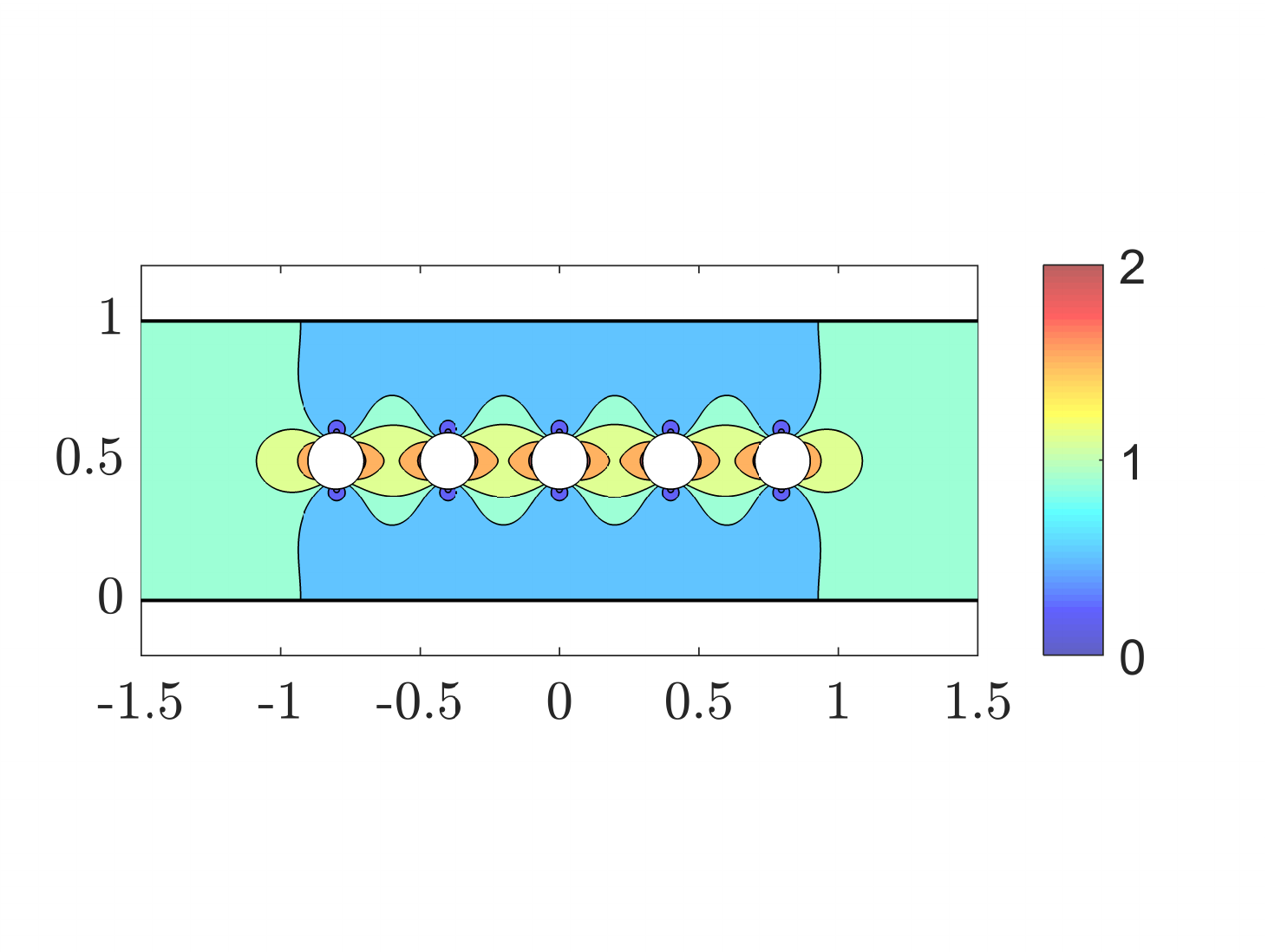}
}
\centerline{
\includegraphics[page=1, trim =0.65cm 2.25cm 0.65cm 2.75cm, clip, width=0.52\textwidth]{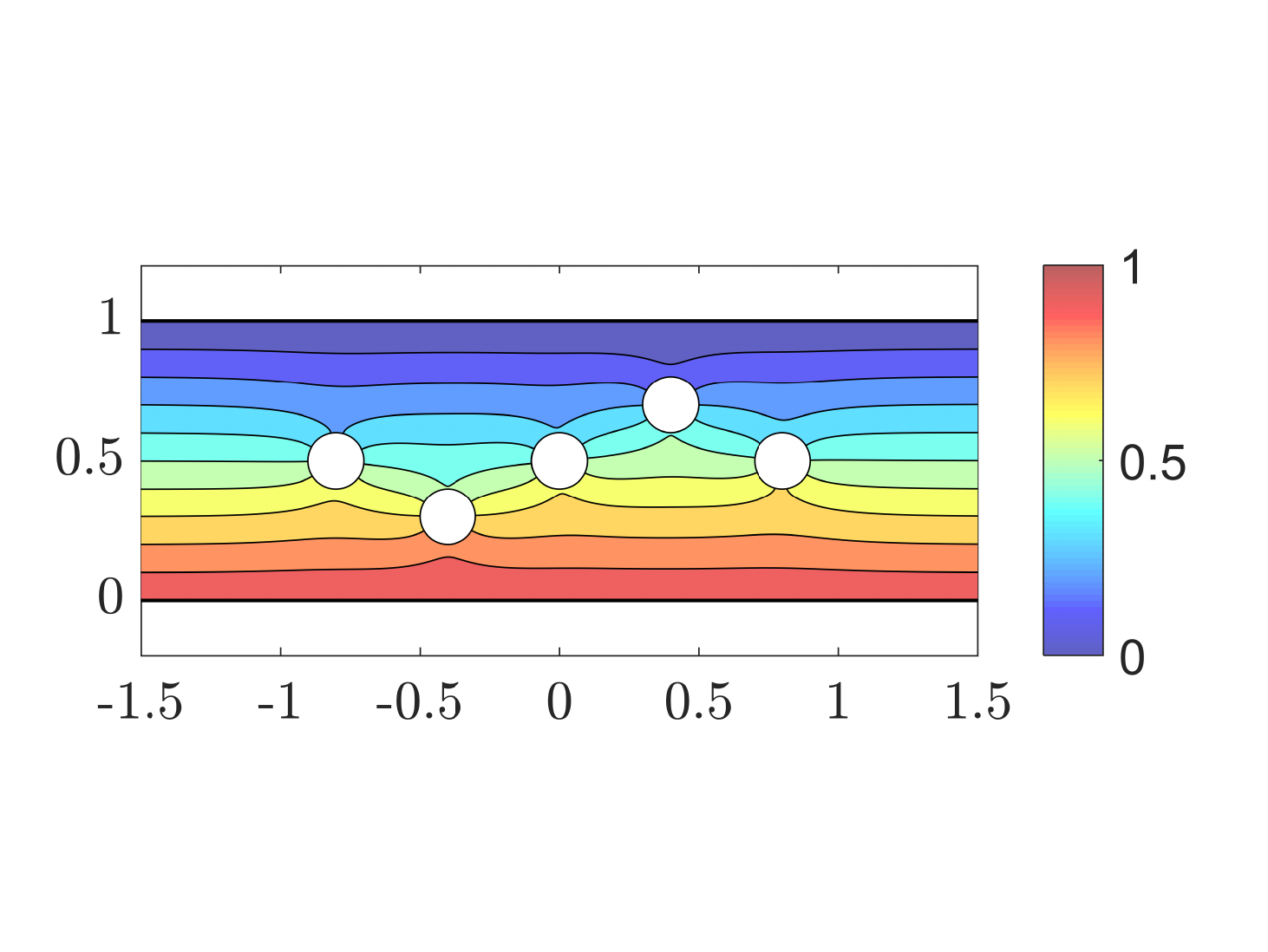}
\hfill
\includegraphics[page=1, trim =0.65cm 2.25cm 0.65cm 2.75cm, clip, width=0.52\textwidth]{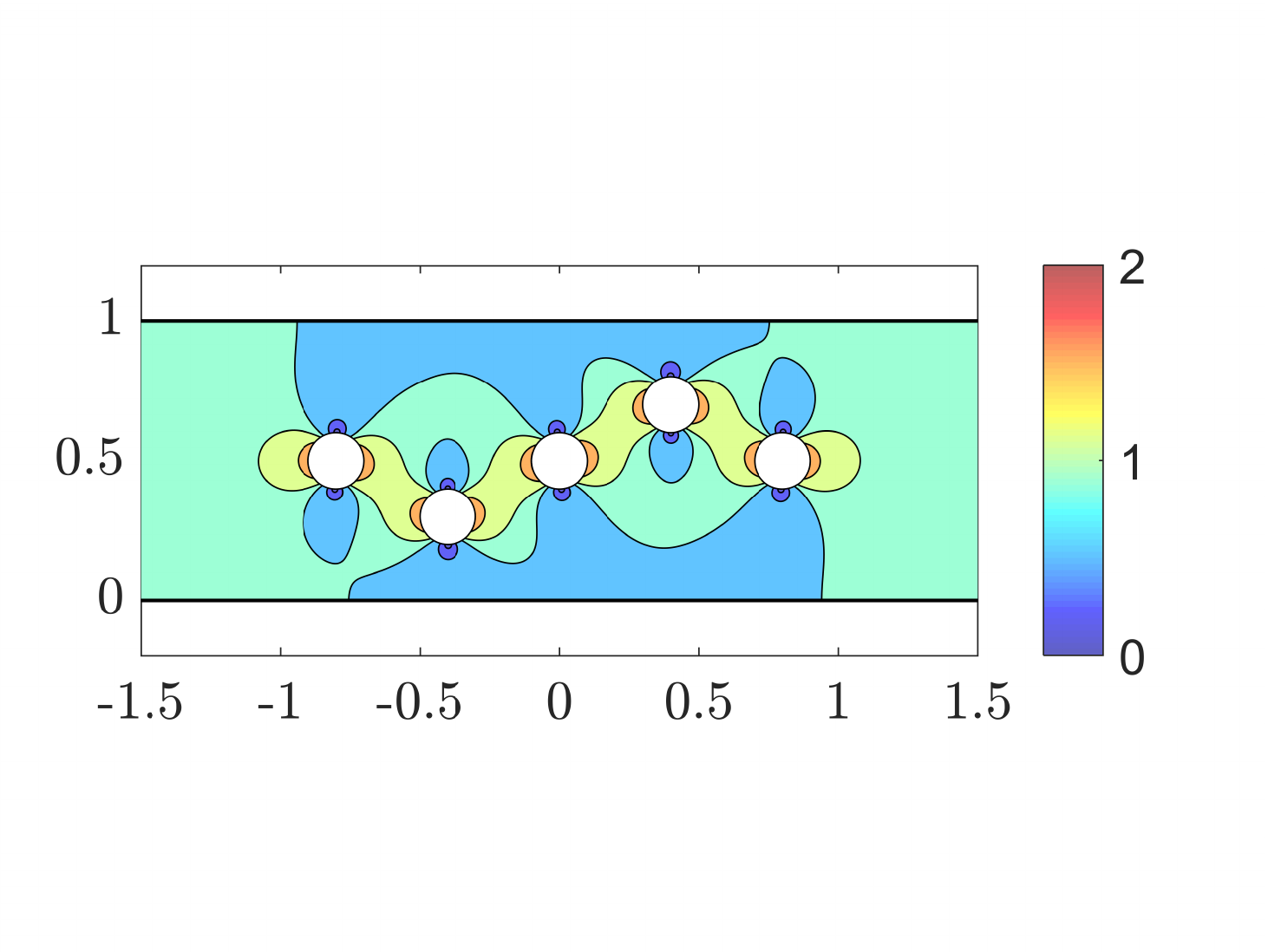}
}
\caption{A contour plot of the temperature distribution $T$ and the heat flux $|q|$ for the domain $\Omega$ with $m=5$ circular holes (Example~\ref{ex:cir-5} for $r=0.1$). First row for Case I and second row for Case II.}
\label{fig:cir-5}
\end{figure}

\begin{figure}[ht] %
\centerline{
\includegraphics[page=1, trim =0cm 0cm 0cm 0cm, clip, width=0.5\textwidth]{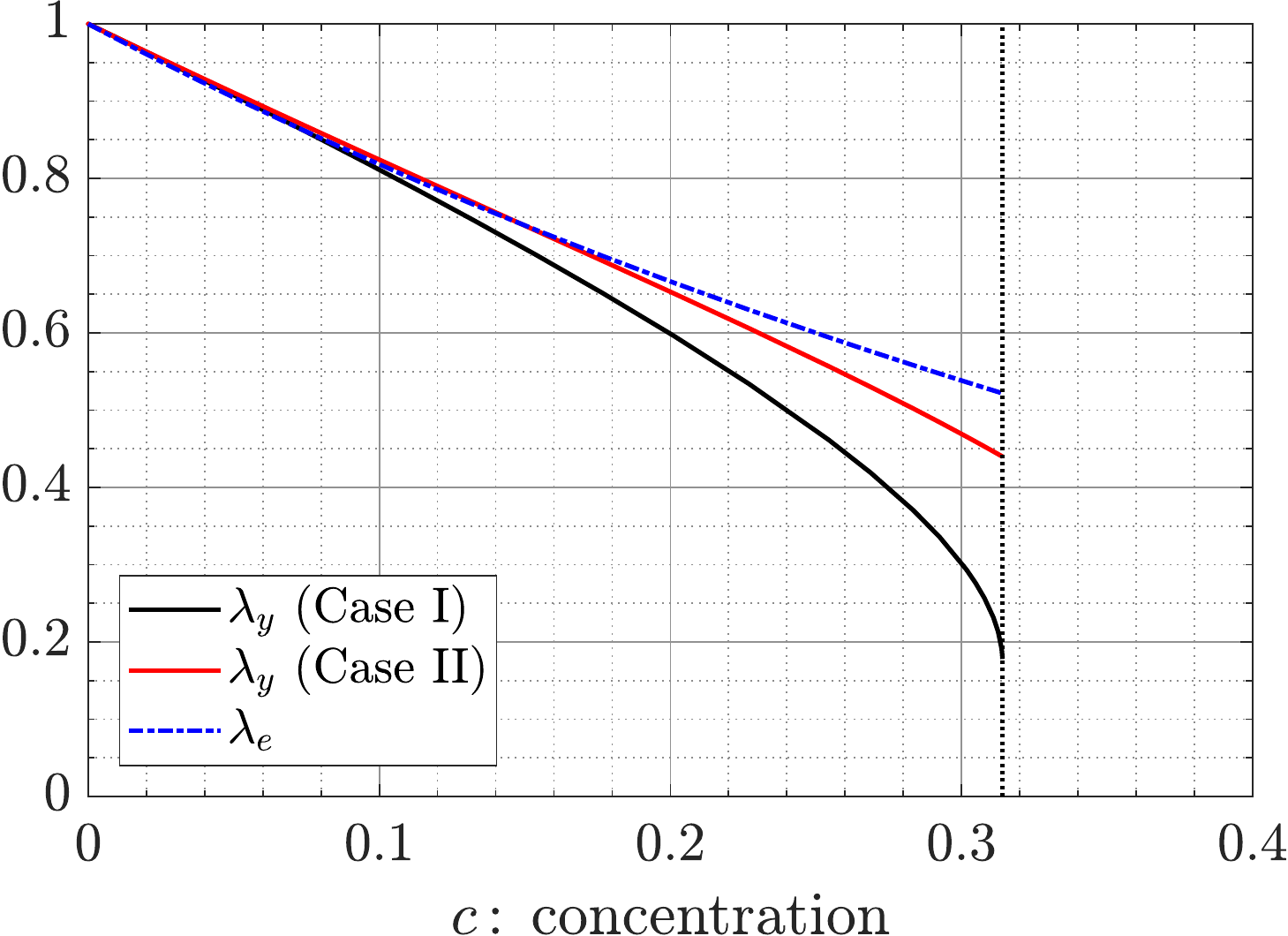}
}
\caption{The effective thermal conductivity $\lambda_y$ and the estimated effective conductivity $\lambda_e$ in~\eqref{eq:CMA} vs. the concentration $c(m,r)=5 r^2\pi/2$ for the domain $\Omega$ with $m=5$ circular holes for $0.00001\le r\le 0.19999$. The vertical dotted line is $c=\pi/10$.}
\label{fig:cir-5L}
\end{figure}

\begin{example}\label{ex:cir-5}{\rm
We consider $m=5$ circular holes with the radius $r$ for $0<r<0.2$. For Case I, we assume the centers of the holes to be set to $-0.8+0.5\i$, $-0.4+0.5\i$, $0.5\i$, $0.4+0.5\i$, and $0.8+0.5\i$.
The contour plot of $T$ and $|q|$ for $r=0.1$ are shown in Figure~\ref{fig:cir-5} (first row).
The approximate value of the effective thermal conductivity for $r=0.1$ is
\[
\lambda_y = 0.8533491.
\]

When $r$ is close to $0.2$, the circular holes become adjacent to each other. To show the effects of the radius $r$ on the effective thermal conductivity $\lambda_y$, we compute the values of $\lambda_y$ for several values of $r$, $0.00001\le r\le 0.19999$. The obtained results are presented in Figure~\ref{fig:cir-5L} where, by~\eqref{eq:cir-conc}, the concentration of these $5$ holes is $c=c(5,r) = 5 r^2\pi/2 \approx 7.854r^2$ for $0<c<\pi/10$ and $0<r<0.2$. The values of the estimated effective conductivity $\lambda_e$ is given also in Figure~\ref{fig:cir-5L}. As one can expect, there is a good agreement between $\lambda_y$ and $\lambda_e$ for small values of $c$. In the same time, the divergence of  $\lambda_y$ and $\lambda_e$ is observed for the concentrations greater than $0.1$.

For Case II, the centers of the holes become $-0.8+0.5\i$, $-0.4+0.3\i$, $0.5\i$, $0.4+0.7\i$, and $0.8+0.5\i$, which means the centers are not anymore horizontally aligned as the second and fourth centers are now shifted by 0.2 up and down, respectively. This is displayed in Figure~\ref{fig:cir-5} (second row).
The curve showing the obtained values of $\lambda_y$ as a function of the concentration is depicted in Figure~\ref{fig:cir-5L}.

Figure~\ref{fig:cir-5L} illustrates that the values of $\lambda_y$ depend on the position of the circular holes centers while the values of $\lambda_e$ are the same for both cases since it depends only on the concentration of the circular holes and not on their positions. We notice a better agreement between $\lambda_y$ and $\lambda_e$ in Cases II when comparing to Case I.
}
\end{example}

\begin{figure}[ht] %
\centerline{
\includegraphics[page=1, trim =0.65cm 2.25cm 0.65cm 2.75cm, clip, width=0.52\textwidth]{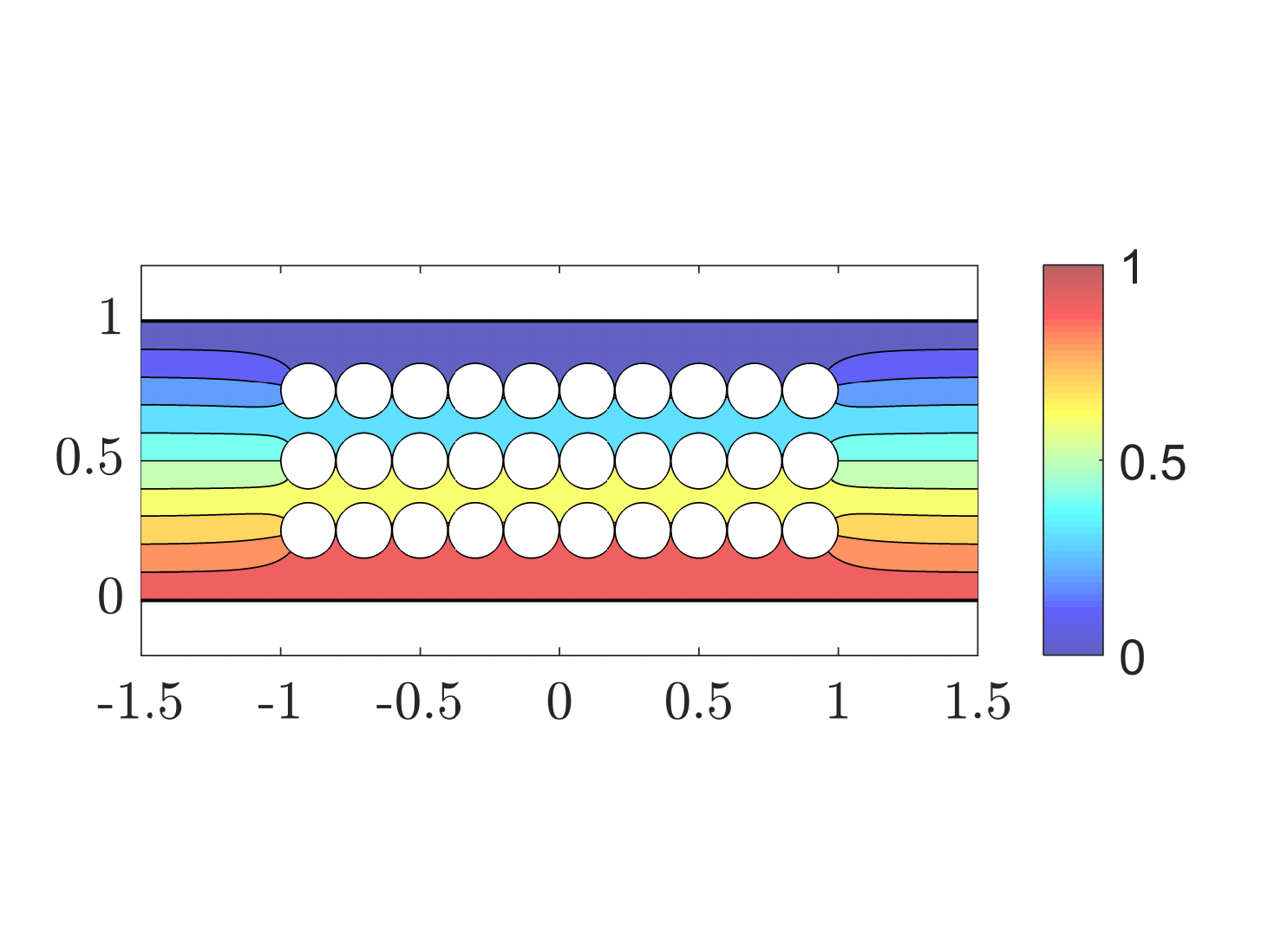}
\hfill
\includegraphics[page=1, trim =0.65cm 2.25cm 0.65cm 2.75cm, clip, width=0.52\textwidth]{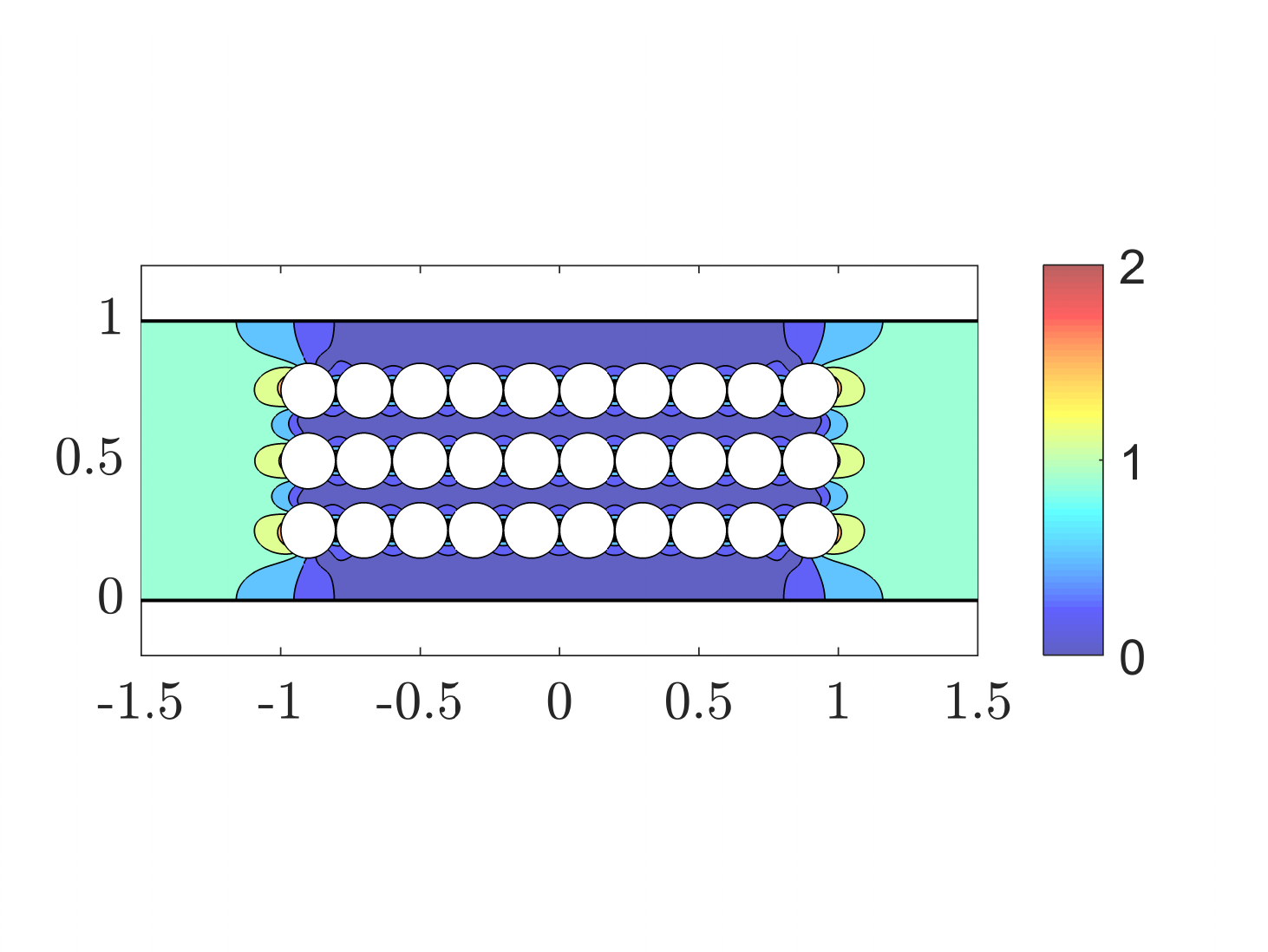}
}
\caption{A contour plot of the temperature distribution $T$ and the heat flux $|q|$ for the domain $\Omega$ with $30$ circular holes ($r=0.099$).}
\label{fig:cir-30}
\end{figure}

\begin{figure}[ht] %
\centerline{
\includegraphics[page=1, trim =0cm 0cm 0cm 0cm, clip, width=0.5\textwidth]{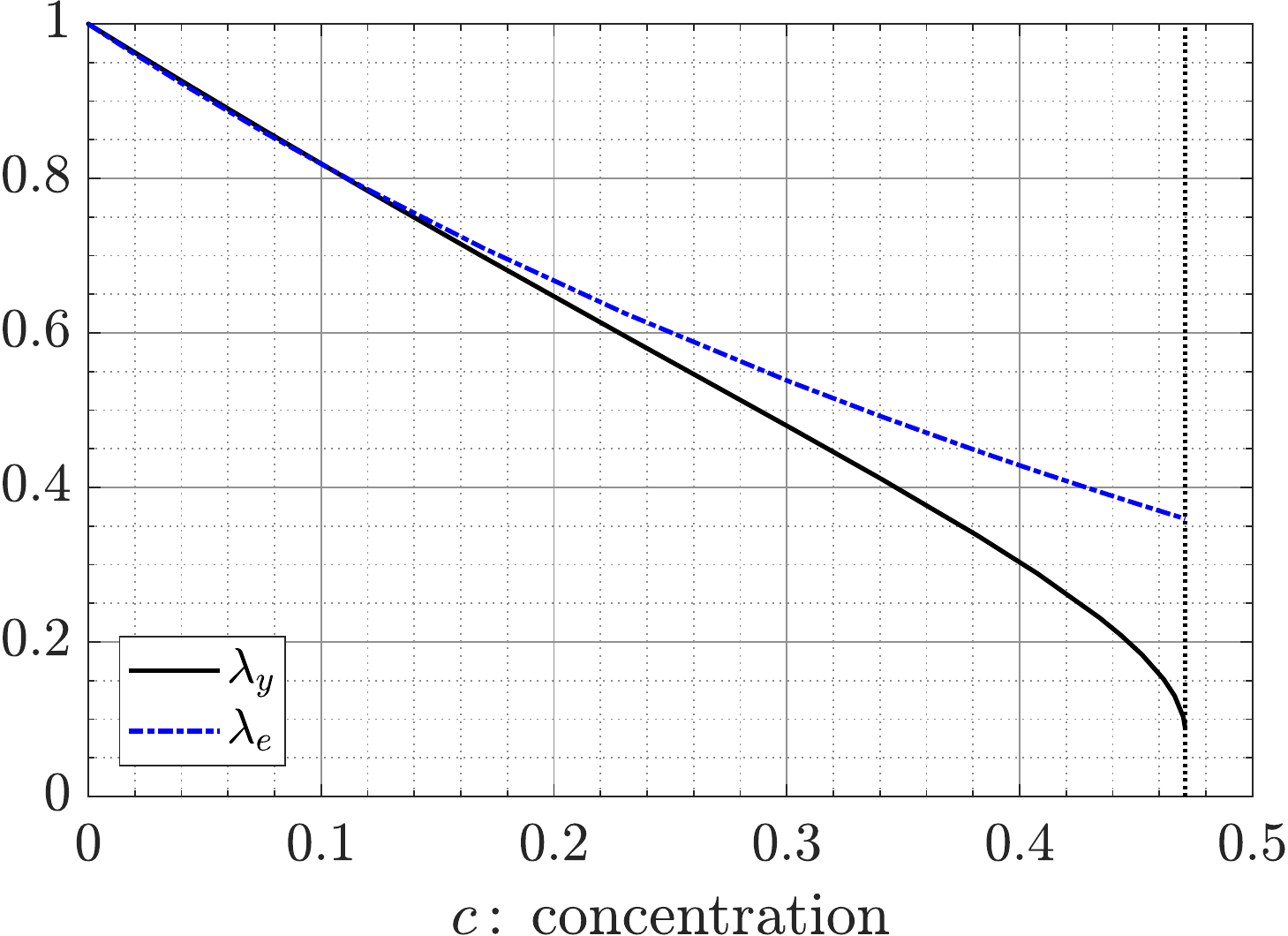}
\hfill
\includegraphics[page=1, trim =0cm 0cm 0cm 0cm, clip, width=0.5\textwidth]{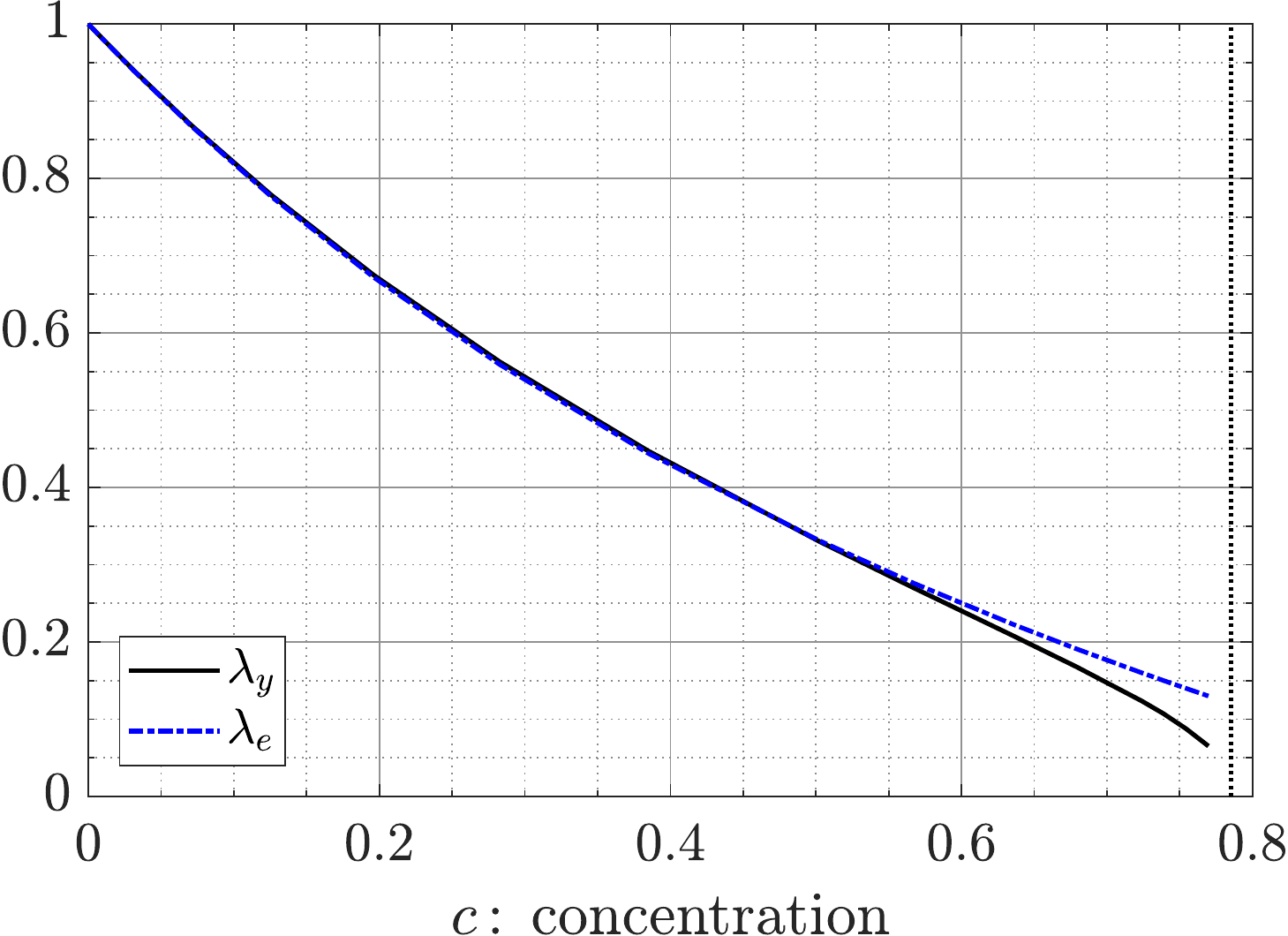}
}
\caption{The effective thermal conductivity $\lambda_y$ vs. the concentration $c(m,r)=m r^2\pi/2$. On the left, the domain $\Omega$ with $m=30$ circular holes (Example~\ref{ex:cir-30}) and $0.00001\le r\le 0.09999$. The vertical dotted line is $c=3\pi/20$. On the right, the domain $\Omega$ with $m=50$ circular holes (Example~\ref{ex:cir-50}) and $0.001\le r\le 0.099$. The vertical dotted line is $c=\pi/4$.}
\label{fig:cir-30L}
\end{figure}

\begin{example}\label{ex:cir-30}{\rm
We consider $m=30$ circular holes with centers $x_k+0.25\i$, $x_k+0.5\i$, and $x_k+0.75\i$, where $x_k=-0.9+0.2(k-1)$ for $k=1,2,\ldots,10$, and with radius $r$ for $0<r<0.1$. 
The contour plot of $T$ and $|q|$ for $r=0.099$ are shown in Figure~\ref{fig:cir-30}.
The approximate value of the effective thermal conductivity for $r=0.099$ is
\[
\lambda_y = 0.1519156.
\]

When $r$ is close to $0.1$, the circular holes become adjacent to each other. We compute the values of $\lambda_y$ for several values of $r$, $0.00001\le r\le 0.09999$. The obtained results are depicted in Figure~\ref{fig:cir-30L} (left) where, by~\eqref{eq:cir-conc}, the concentration of these $30$ holes is $c=c(30,r) = 30 r^2\pi/2\approx 47.124r^2$. Note that $0<c<3\pi/20$ for $0<r<0.1$.
}\end{example}

\begin{example}\label{ex:cir-50}{\rm
We take up here the case of $m=50$ circular holes with centers $x_k+0.1\i$, $x_k+0.3\i$, $x_k+0.5\i$, $x_k+0.7\i$, and $x_k+0.9\i$, where $x_k=-0.9+0.2(k-1)$ for $k=1,2,\ldots,10$, and with radius $r$ for $0<r<0.1$. 
On the basis of~\eqref{eq:cir-conc}, the concentration of these $50$ holes is $c=c(50,r) = 50 r^2\pi/2\approx 78.54r^2$. For $0<r<0.1$, we have  $0<c<\pi/4$. When $r$ is close to $0.1$, the circular holes become adjacent to each other, and the concentration is almost equal to $\pi/4$. The obtained results showing the behavior of $\lambda_y$ as a function of the radius $r$, for $0.001\le r\le 0.099$, are presented in Figure~\ref{fig:cir-30L} (right). 
}\end{example}

\subsection{The domain $\Omega$ with only CNTs}

In this subsection, we consider the domain $\Omega$ with $m$ non-overlapping elliptic CNTs without any circular holes (i.e., $m=\ell$). We assume that all CNTs have equal sizes and are of elliptic shape where the ellipses have the parametrization
\begin{equation}\label{eq:ellipse}
\eta_j(t) = z_j+ a\cos t-\i b \sin t, \quad 0\le t\le 2\pi, \quad j=1,2,\ldots,m,
\end{equation}
where $z_j$ is the center of the ellipse, $2a$ and $2b$ are the length of the ellipses axes in the $x$ and $y$-directions, respectively. If $a/b>1$, the major axis of the ellipses is horizontal, if $a/b<1$, the major axis of the ellipses is vertical, and if $a/b=1$,  the ellipses reduced to circles. Here, we choose $a$ and $b$ such that their ratio satisfies $0.1\le a/b\le 10$.
These elliptic shape CNTs are chosen in the part of the domain $\Omega$ between $x=-1$ and $x=1$. So, we define the concentration $c(m,a,b)$ of these CNTs to be 
\begin{equation}\label{eq:ell-conc}
c(m,a,b) = \frac{m ab \pi}{2}.
\end{equation}
If $\frac{a}{b} \ll 1$, instead of \eqref{eq:ell-conc} the plane slits density is considered in the theory of composites and porous media
\begin{equation}\label{eq:ell-concP}
\phi = \frac{m b^2}{|\Omega|}= \frac{m b^2}{2},
\end{equation}

	For a macroscopically isotropic media with only perfectly conducting identical circular inclusions (CNTs), an approximation of the effective conductivity $\lambda_e$ is given by the inverse to \eqref{eq:CMA} value (see \cite{MR})
	\begin{equation}\label{eq:CMA2}
		\lambda_e = \frac{1+c}{1-c} + O(c^3).
	\end{equation}

\begin{figure}[ht] %
\centerline{
\includegraphics[page=1, trim =0.65cm 2.25cm 0.65cm 2.75cm, clip, width=0.52\textwidth]{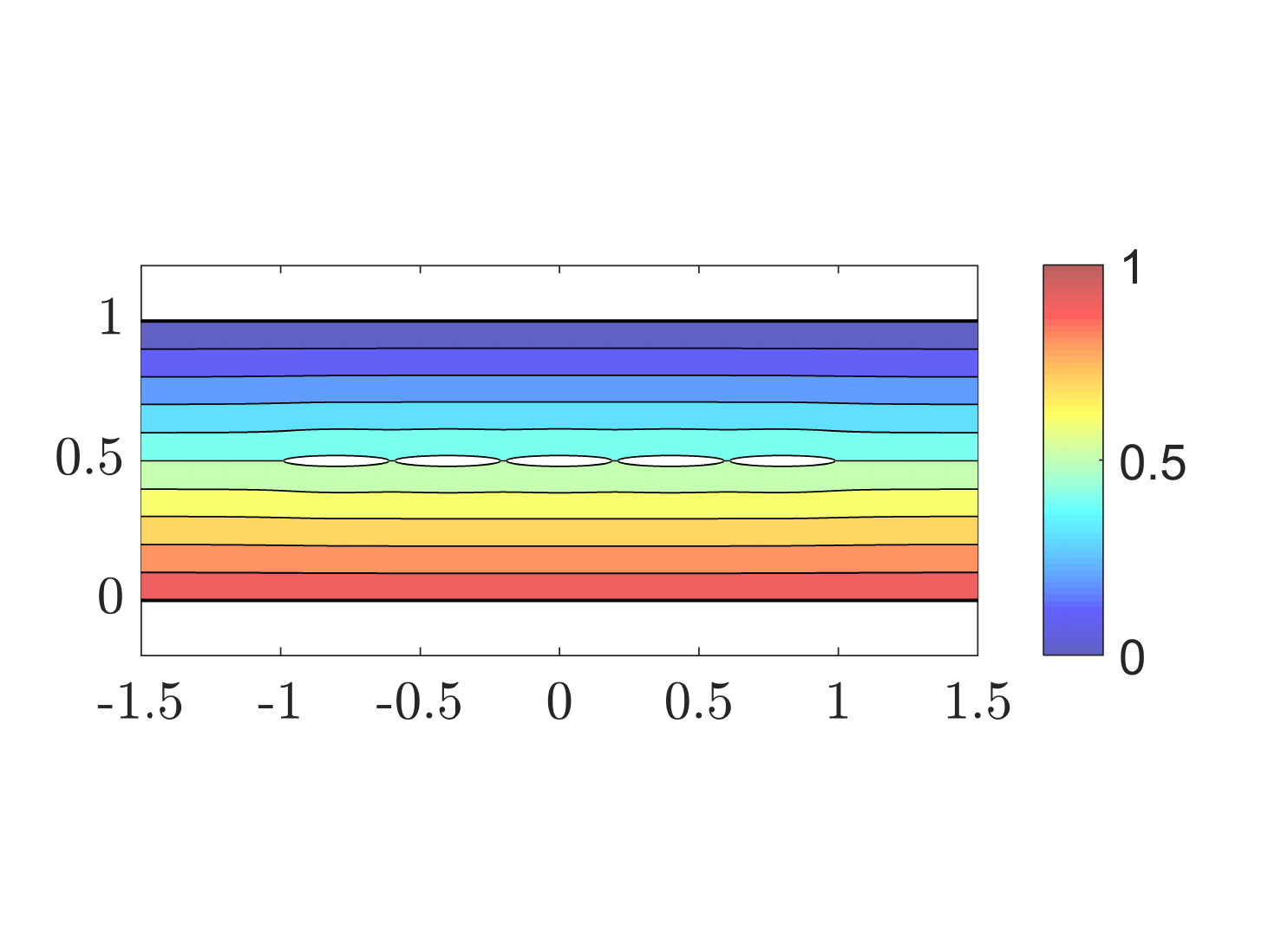}
\hfill
\includegraphics[page=1, trim =0.65cm 2.25cm 0.65cm 2.75cm, clip, width=0.52\textwidth]{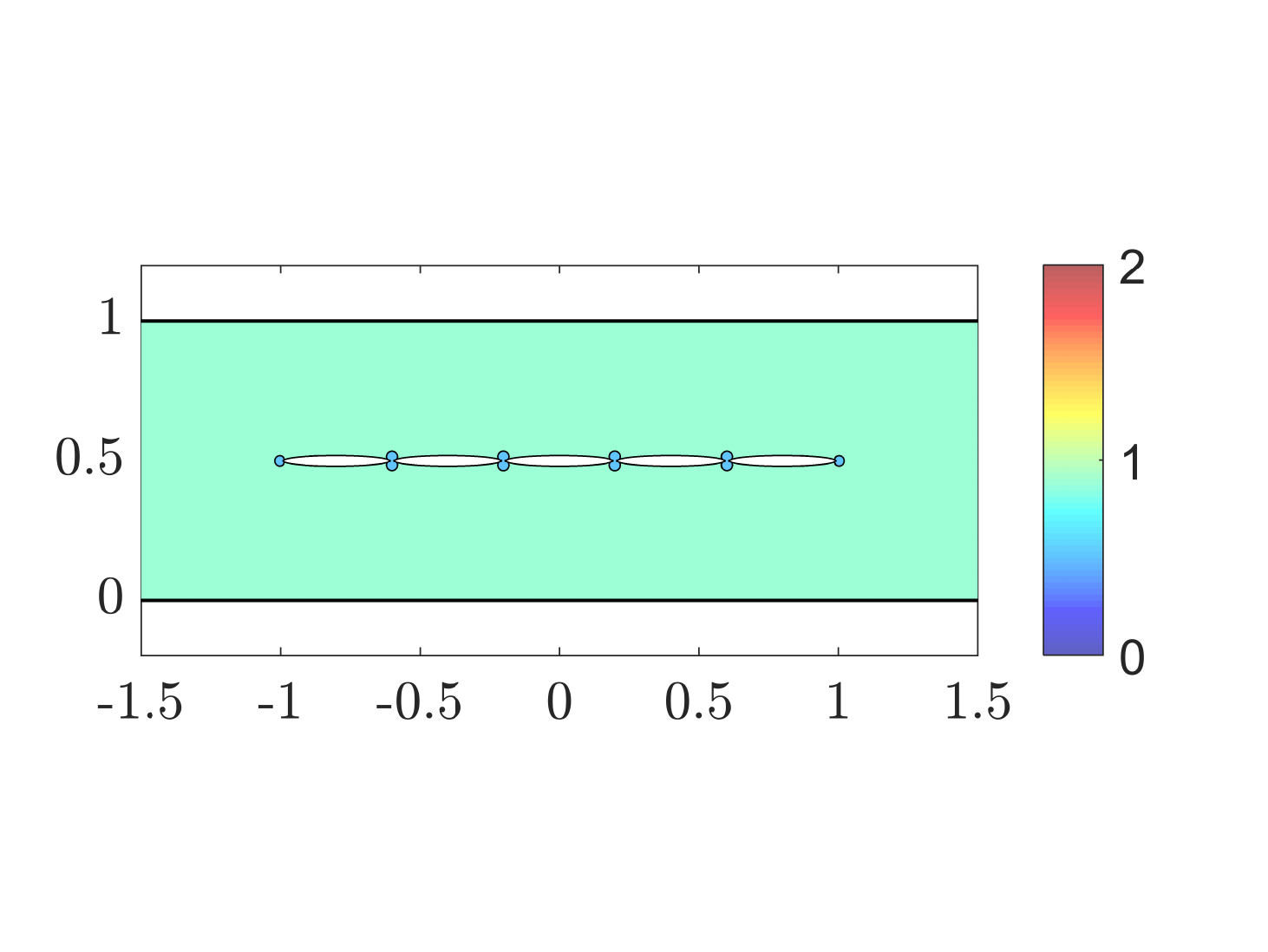}
}
\centerline{
\includegraphics[page=1, trim =0.65cm 2.25cm 0.65cm 2.75cm, clip, width=0.52\textwidth]{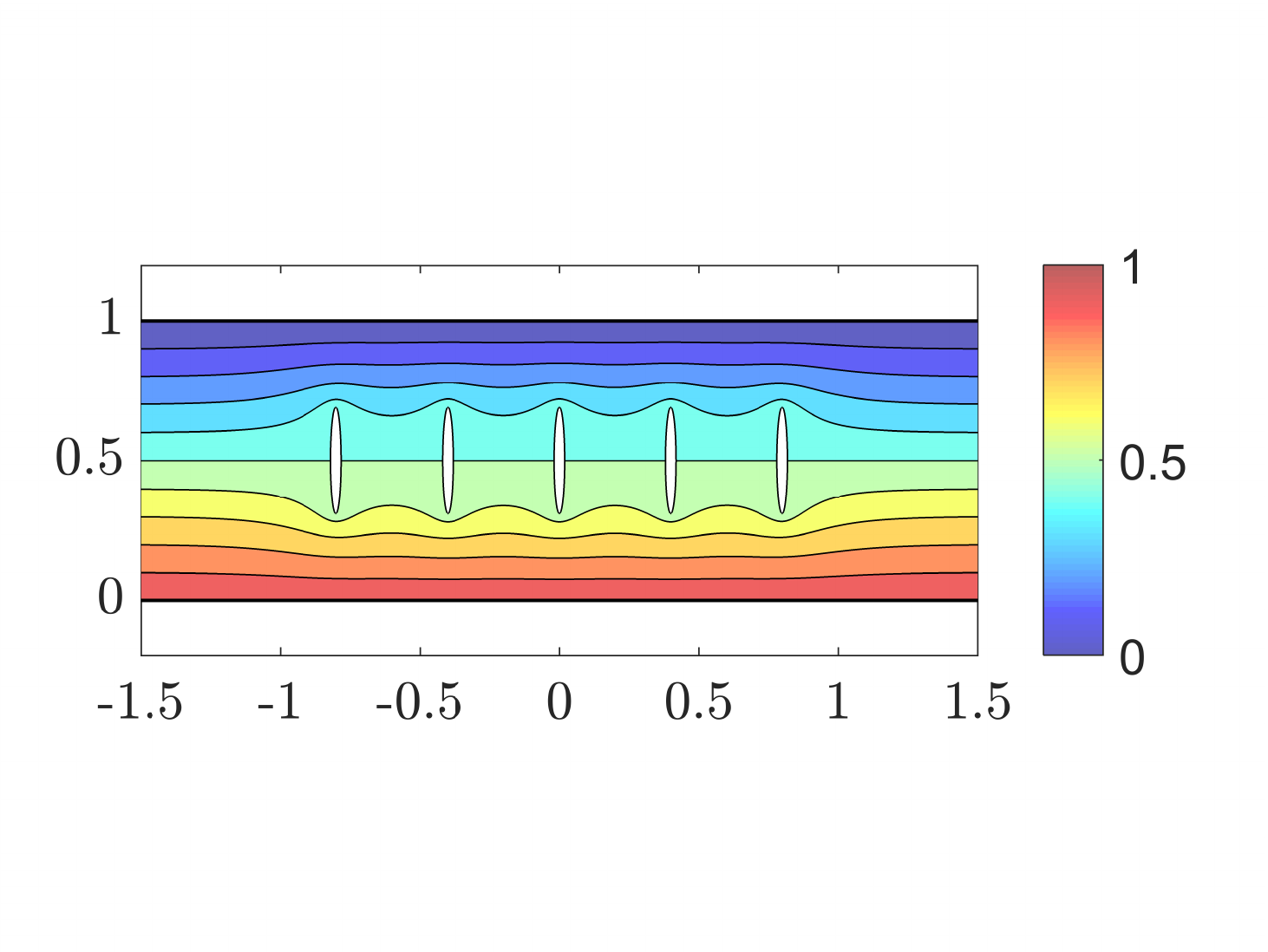}
\hfill
\includegraphics[page=1, trim =0.65cm 2.25cm 0.65cm 2.75cm, clip, width=0.52\textwidth]{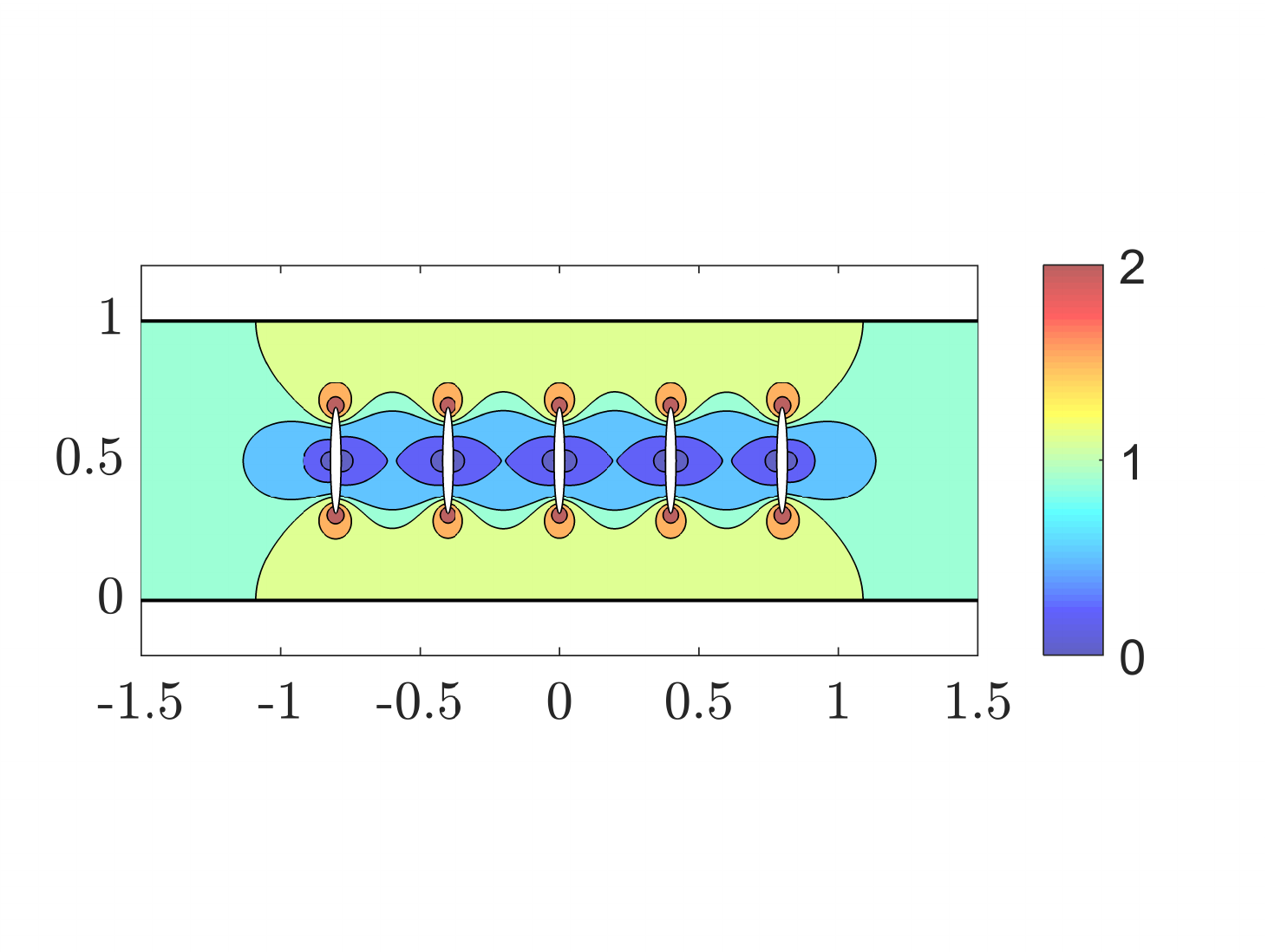}
}
\caption{A contour plot of the temperature distribution $T$ and the heat flux $|q|$ for the domain $\Omega$ with $m=5$ elliptic CNTs (Example~\ref{ex:ell-5}), where $a=0.19$ and $b=0.019$ for the first row and $a=0.019$ and $b=0.19$ for the second row.}
\label{fig:ell-5e}
\end{figure}

\begin{figure}[ht] %
\centerline{
\includegraphics[page=1, trim =0cm 0cm 0cm 0cm, clip, width=0.5\textwidth]{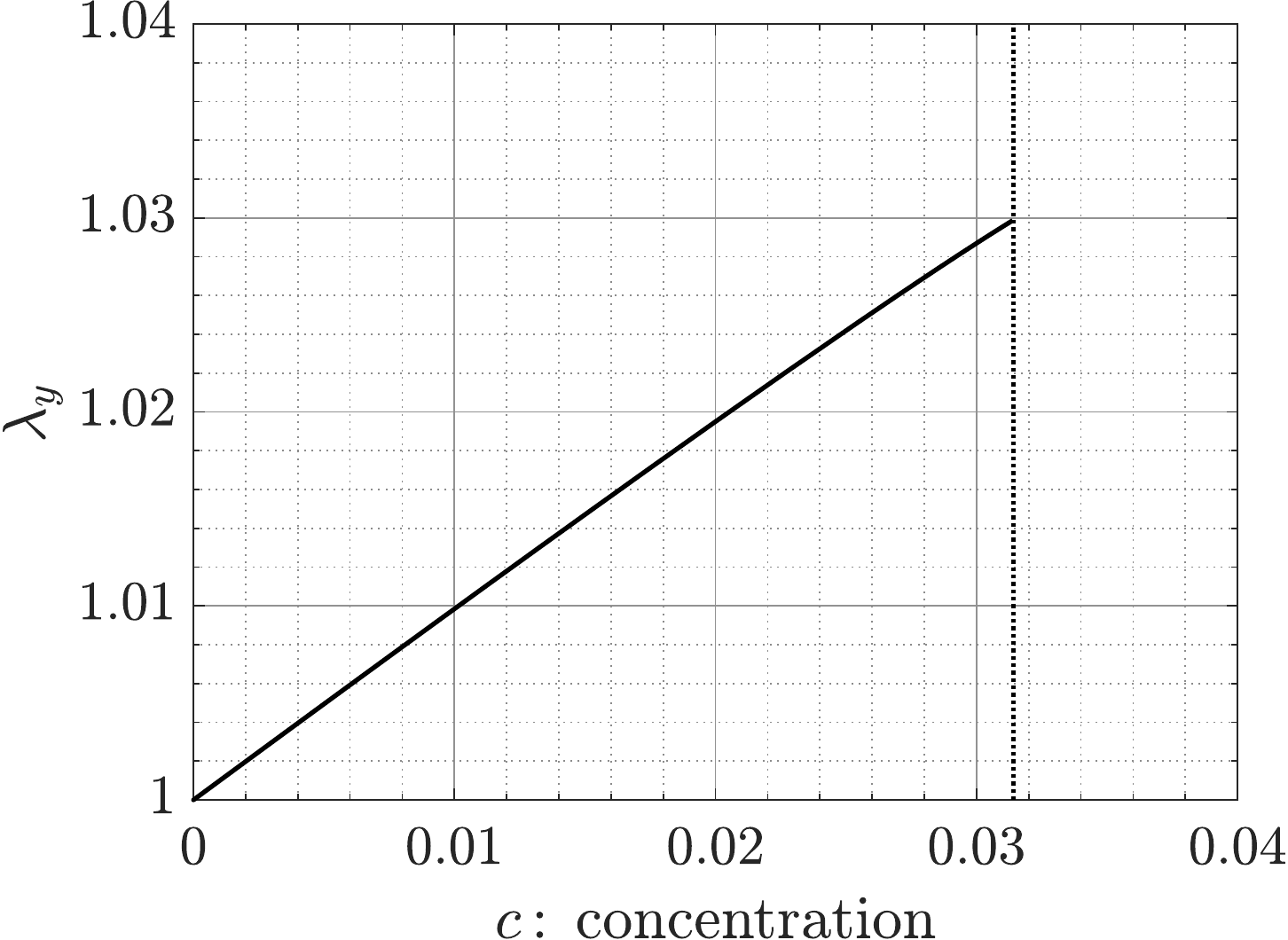}
\hfill
\includegraphics[page=1, trim =0cm 0cm 0cm 0cm, clip, width=0.5\textwidth]{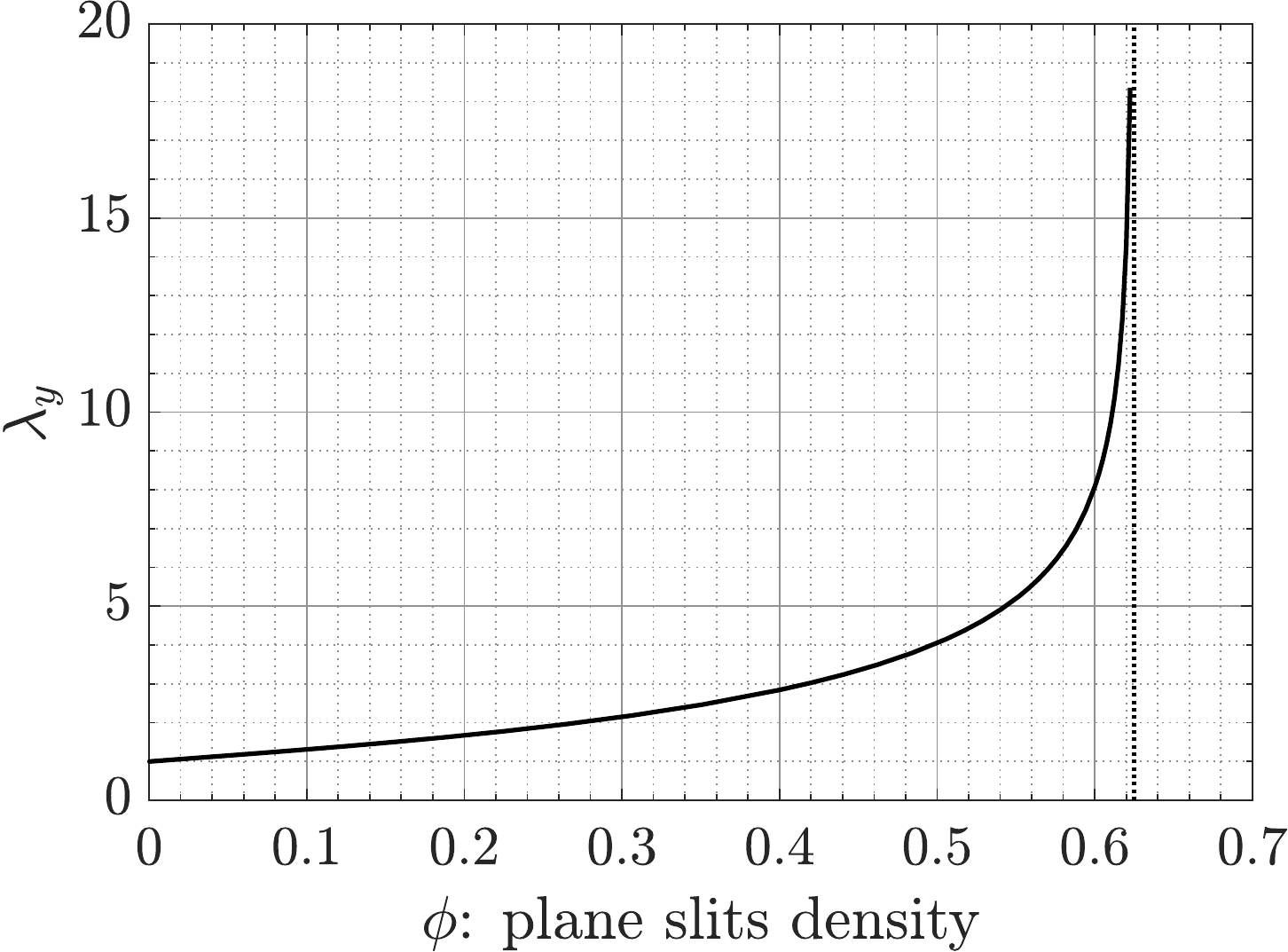}
}
\caption{The effective thermal conductivity $\lambda_y$ for the domain $\Omega$ with $m=5$ elliptic CNTs (Example~\ref{ex:ell-5}).
On the left, the effective thermal conductivity $\lambda_y$ vs. the concentration $c = a^2\pi/4$ for $0.0001\le a\le 0.1999$ and $a/b=10$. The vertical dotted line is $c=\pi/100$. On the right, the effective thermal conductivity $\lambda_y$ vs. the plane slits density $\phi = 2.5 b^2$ for $0.001\le b\le 0.499$ and $a/b=0.1$. The vertical dotted line is $\phi=0.625$.}
\label{fig:ell-5hv}
\end{figure}

\begin{example}\label{ex:ell-5}{\rm
We consider $m=5$ elliptic CNTs with centers $-0.8+0.5\i$, $-0.4+0.5\i$, $0.5\i$, $0.4+0.5\i$, and $0.8+0.5\i$, and where $0<a<0.2$ and $0<b<0.5$. 
Figure~\ref{fig:ell-5e} (first row) presents the contour plot of $T$ and $|q|$ for $a=0.19$ and $b=0.019$ (the ellipses are horizontal). For these values of $a$ and $b$, the approximate value of the effective thermal conductivity is
\[
\lambda_y = 1.0272480.
\]
For $a=0.019$ and $b=0.19$, the contour plot of $T$ and $|q|$ are shown in Figure~\ref{fig:ell-5e} (second row). The approximate value of the effective thermal conductivity for these values of $a$ and $b$ is
\[
\lambda_y = 1.2804116.
\]
The CNTs in Figure~\ref{fig:ell-5e} have the same concentration. However, the value of  $\lambda_y$ is larger for the vertical ellipses case.  

When $a$ approaches $0.2$, the ellipses get adjacent to each other. On the other hand, they come close to the upper and lower walls when $b$ approaches $0.5$.
We compute the values of $\lambda_y$ for several values of $a$, $0.0001\le a\le 0.1999$, with $b=0.1a$. The obtained results are presented in Figure~\ref{fig:ell-5hv} (left) where, by~\eqref{eq:ell-conc}, the concentration of these $5$ ellipses is $c=c(5,a,b) = 5ab\pi/2 = a^2\pi/4\approx 0.7854a^2$. Note that, for $0<a<0.2$ and $b=0.1a$, we have $0<c<\pi/100$.

The values of $\lambda_y$ are also computed for several values of $b$ for $0.001\le b\le 0.499$ with $a=0.1b$. Since $a/b=0.1$ is small, the obtained values of $\lambda_y$ are plotted versus the values of $\phi = 2.5 b^2$, given by~\eqref{eq:ell-concP}, where $0<\phi<5/8$ for $0< b <0.5$. The obtained results are presented in Figure~\ref{fig:ell-5hv} (right).
}\end{example}

\begin{figure}[!ht] %
\centerline{
\includegraphics[page=1, trim =0cm 0cm 0cm 0cm, clip, width=0.524\textwidth]{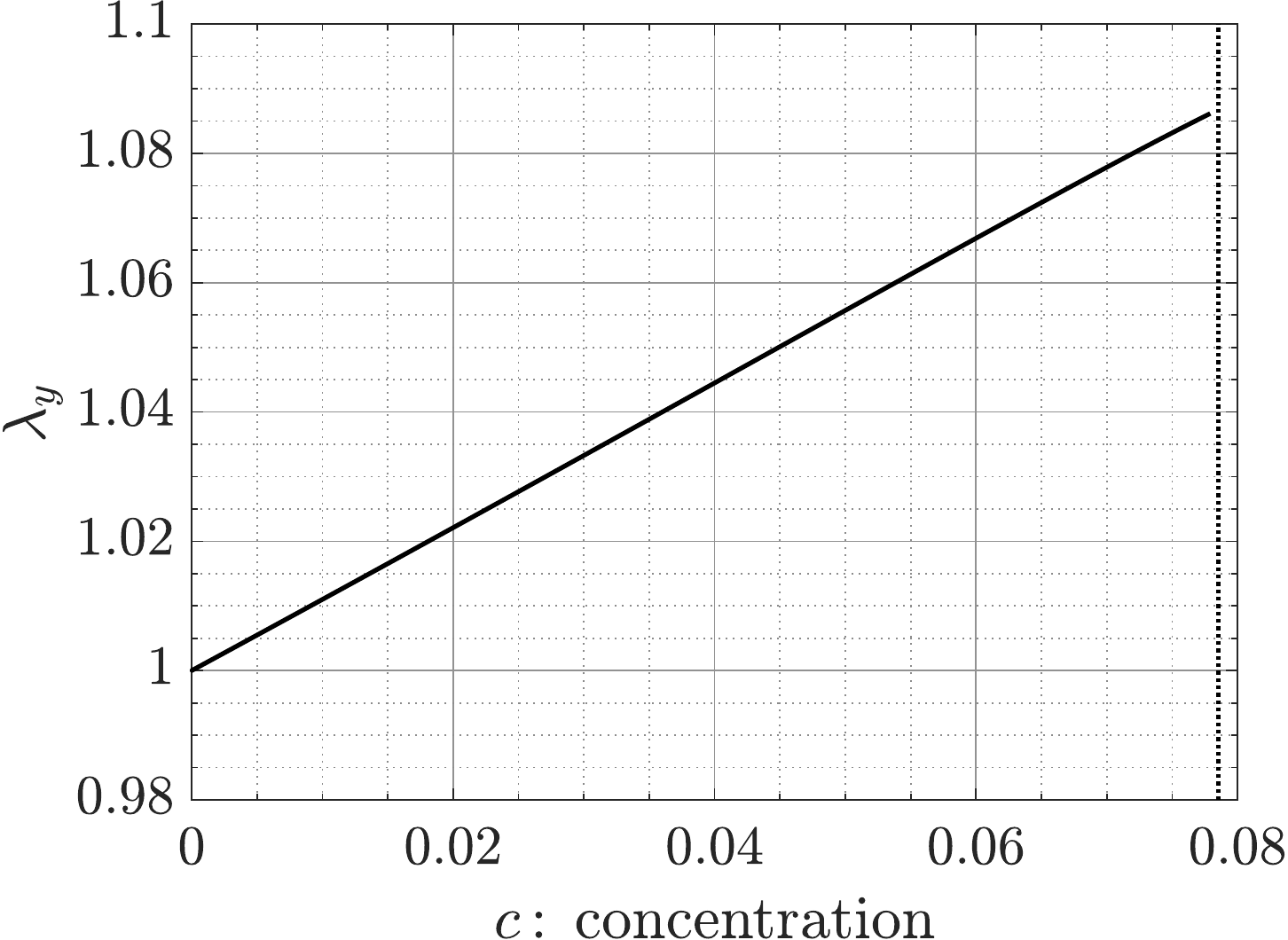}
\hfill
\includegraphics[page=1, trim =0cm 0cm 0cm 0cm, clip, width=0.5\textwidth]{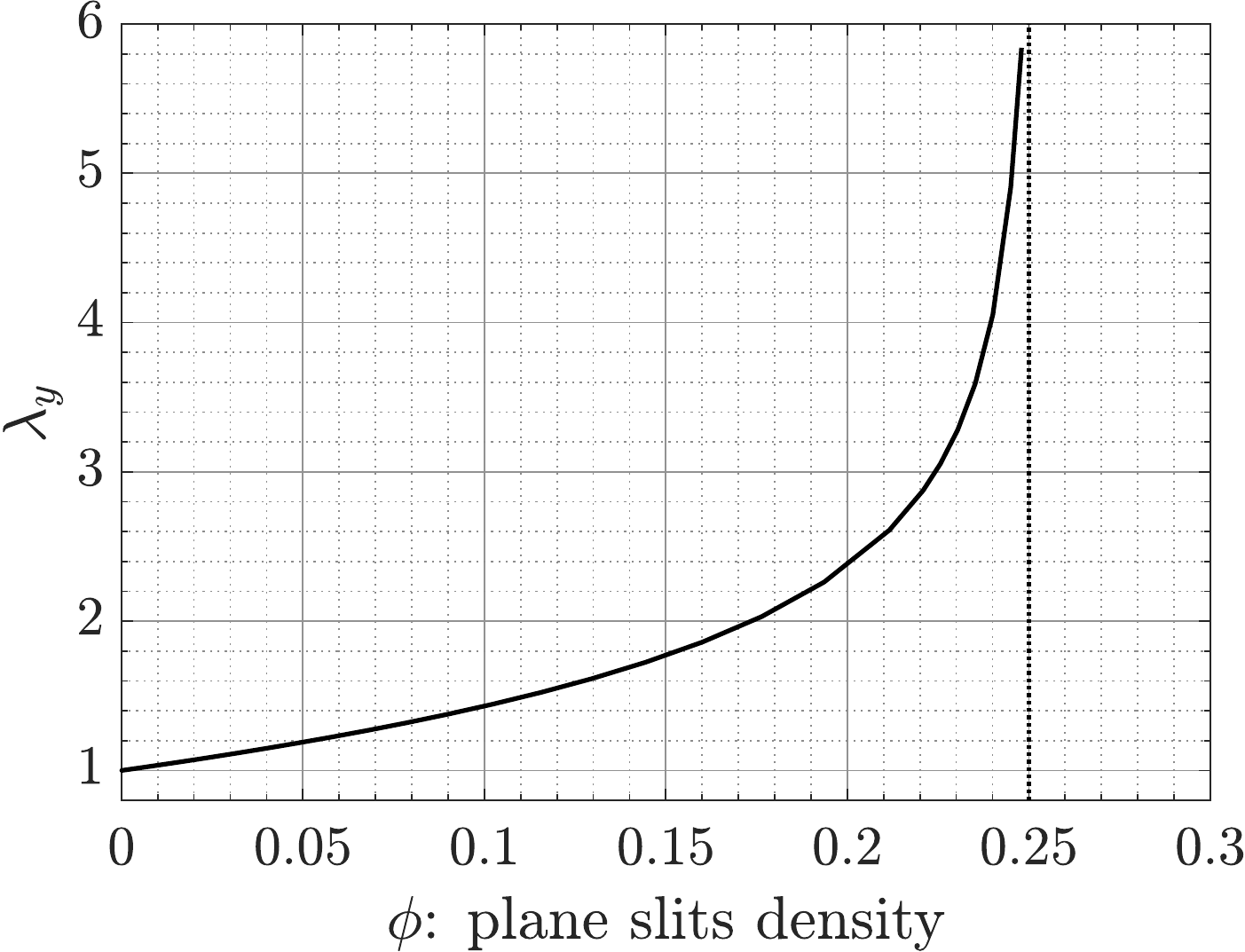}
}
\caption{The effective thermal conductivity $\lambda_y$ for the domain $\Omega$ with $m=200$ elliptic CNTs (Example~\ref{ex:ell-200}).
On the left, the effective thermal conductivity $\lambda_y$ vs. the concentration $c=10a^2\pi$ for $0.0002\le a\le 0.0498$ with $a/b=10$. The vertical dotted line is $c=\pi/40$. On the right, the effective thermal conductivity $\lambda_y$ vs. the plane slits density $\phi= 100 b^2$ for $0.0002\le b\le 0.0498$, with $a/b=0.1$. The vertical dotted line is $\phi=1/4$. The vertical dotted line is $c=1/4$.}
\label{fig:ell-200hv}
\end{figure}

\begin{figure}[ht] %
\centerline{
\includegraphics[page=1, trim =0cm 0cm 0cm 0cm, clip, width=0.5\textwidth]{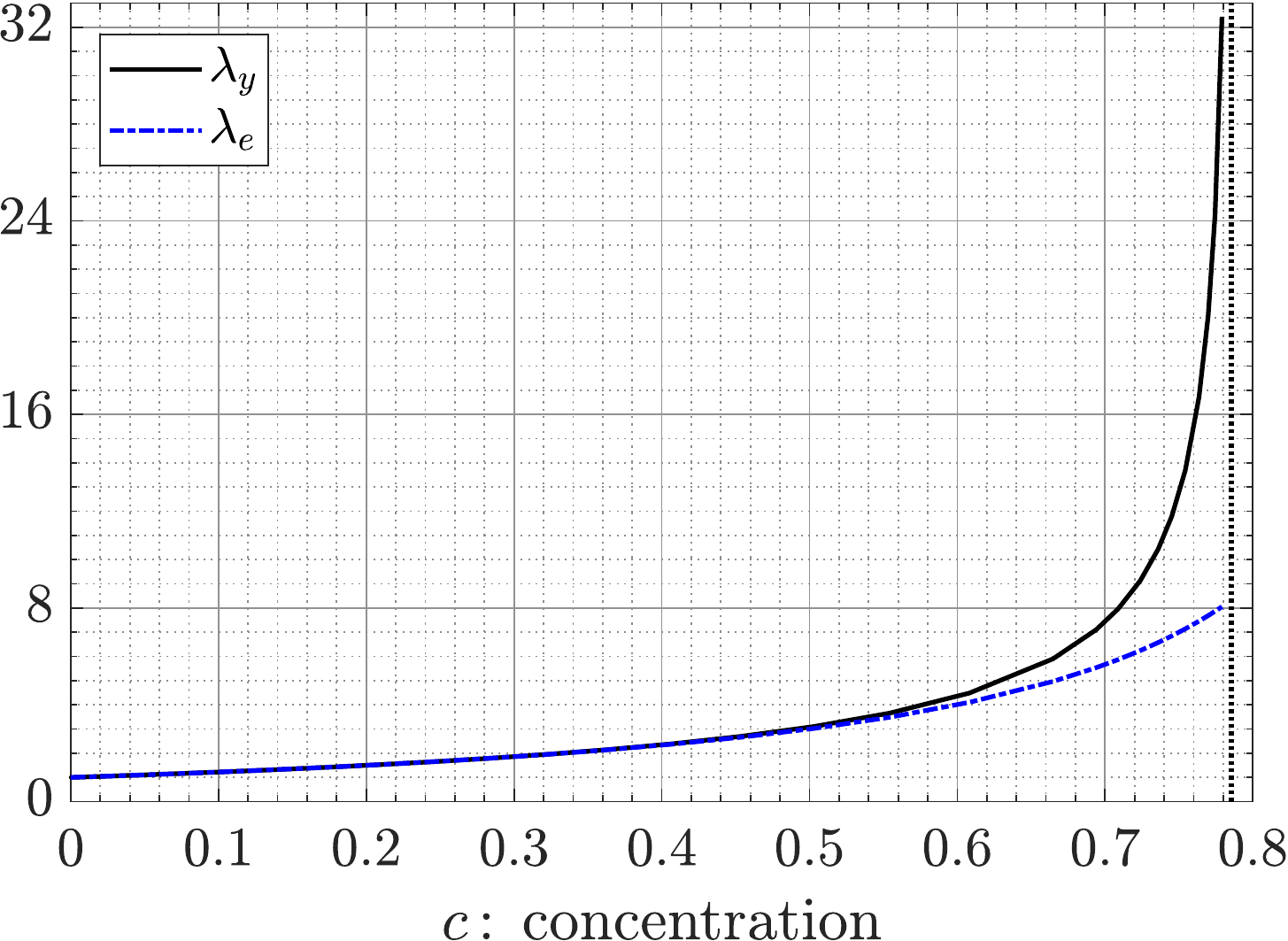}
}
\caption{
The effective thermal conductivity $\lambda_y$ (for the domain $\Omega$ with $m=200$ circular CNTs obtained by setting $b=a$ in Example~\ref{ex:ell-200}) and the estimated effective conductivity $\lambda_e$ in~\eqref{eq:CMA2} vs. the concentration $c=100a^2\pi$ for $0.0002\le a\le 0.0498$. The vertical dotted line is $c=\pi/4$.}
\label{fig:ell-200c}
\end{figure}

\begin{example}\label{ex:ell-200}{\rm
We consider $m=200$ elliptic CNTs with centers $x_k+\i y_j$ for $k=1,2,\ldots,20$ and $j=1,2,\ldots,10$ where $x_k=-0.95+(k-1)/10$ and $y_j=0.05+(j-1)/10$, and with $0<a<0.05$ and $0<b<0.05$. 

We compute the values of $\lambda_y$ for several values of $a$, $0.0002\le a\le 0.0498$, and $a/b=10$ (i.e., the ellipses are horizontal) where the ellipses become close to each other when $a$ approaches $0.05$. The obtained results are presented in Figure~\ref{fig:ell-200hv} (left) where the concentration of these $200$ ellipses is $c= 10a^2\pi\approx 31.416a^2$. For $0<a<0.05$ and $a/b=10$, we have $0<c<\pi/40$.
Then, we compute the values of $\lambda_y$ for several values of $b$, $0.0002\le b\le 0.0498$, and $a/b=0.1$. The obtained results for $\lambda_y$ versus the the plane slits density $\phi=m b^2/2=100b^2$ are presented in Figure~\ref{fig:ell-200hv} (right) where $0<\phi<1/4$ for $0<b<0.05$ and $a/b=0.1$.

When $a/b=1$, the ellipses reduce to circles. We compute the values of $\lambda_y$ for several values of $a$, $0.0002\le a\le 0.0498$. The obtained results are presented in Figure~\ref{fig:ell-200c} where the concentration of these $200$ ellipses is $c= 100a^2\pi\approx 314.16a^2$. For $0<a<0.05$ and $a/b=1$, we have $0<c<\pi/4$. Figure~\ref{fig:ell-200c} presents also the values of the estimated effective conductivity $\lambda_e$ given by~\eqref{eq:CMA2}.
}\end{example}

\subsection{The domain $\Omega$ with $2000$ CNTs and/or circular voids}

We are concerned in this section with the study of a large number of perfect conductors and/or insulators. We consider two example where in the first both perfect conductors and insulators have the same circular shape, while in the second, conductors have an elliptic shape and insulators have a circular shape.
The present investigation is useful when studying the impact of geometric shapes on the macroscopic properties of three-phases high contrast media. 

\begin{figure}[ht] %
\centerline{
\includegraphics[page=1, trim=0cm 2cm 0cm 1.5cm, clip, width=0.65\textwidth]{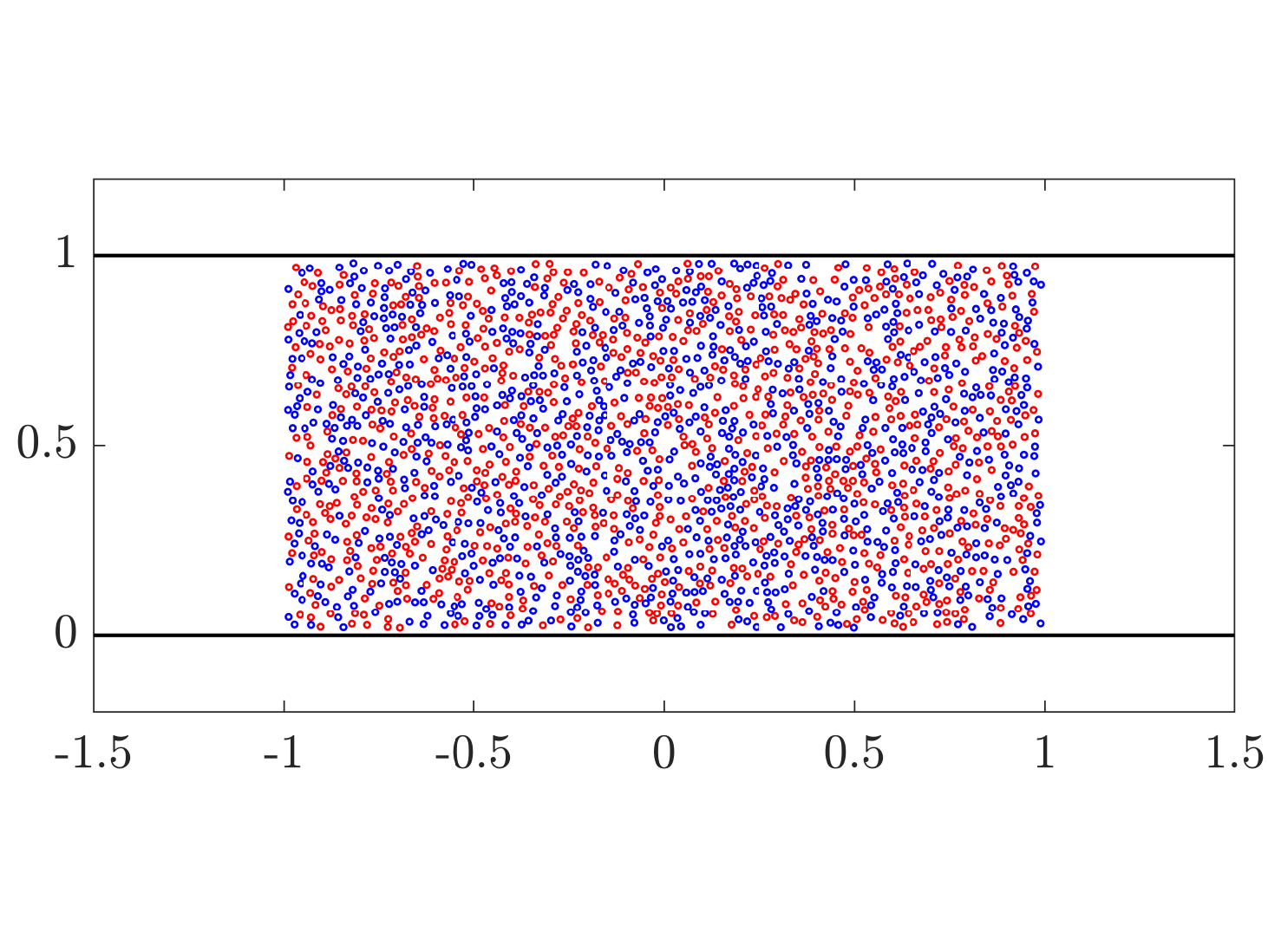}
}
\caption{The domain $\Omega$ with $m=2000$ circular holes. Case I: We have $p=1000$ voids (blue circles) and $\ell=1000$ CNTs (red circles).}
\label{fig:2000NV}
\end{figure}

\begin{figure}[ht] %
\centerline{
\includegraphics[page=1, trim=0cm 0cm 0cm 0cm, clip, width=0.5\textwidth]{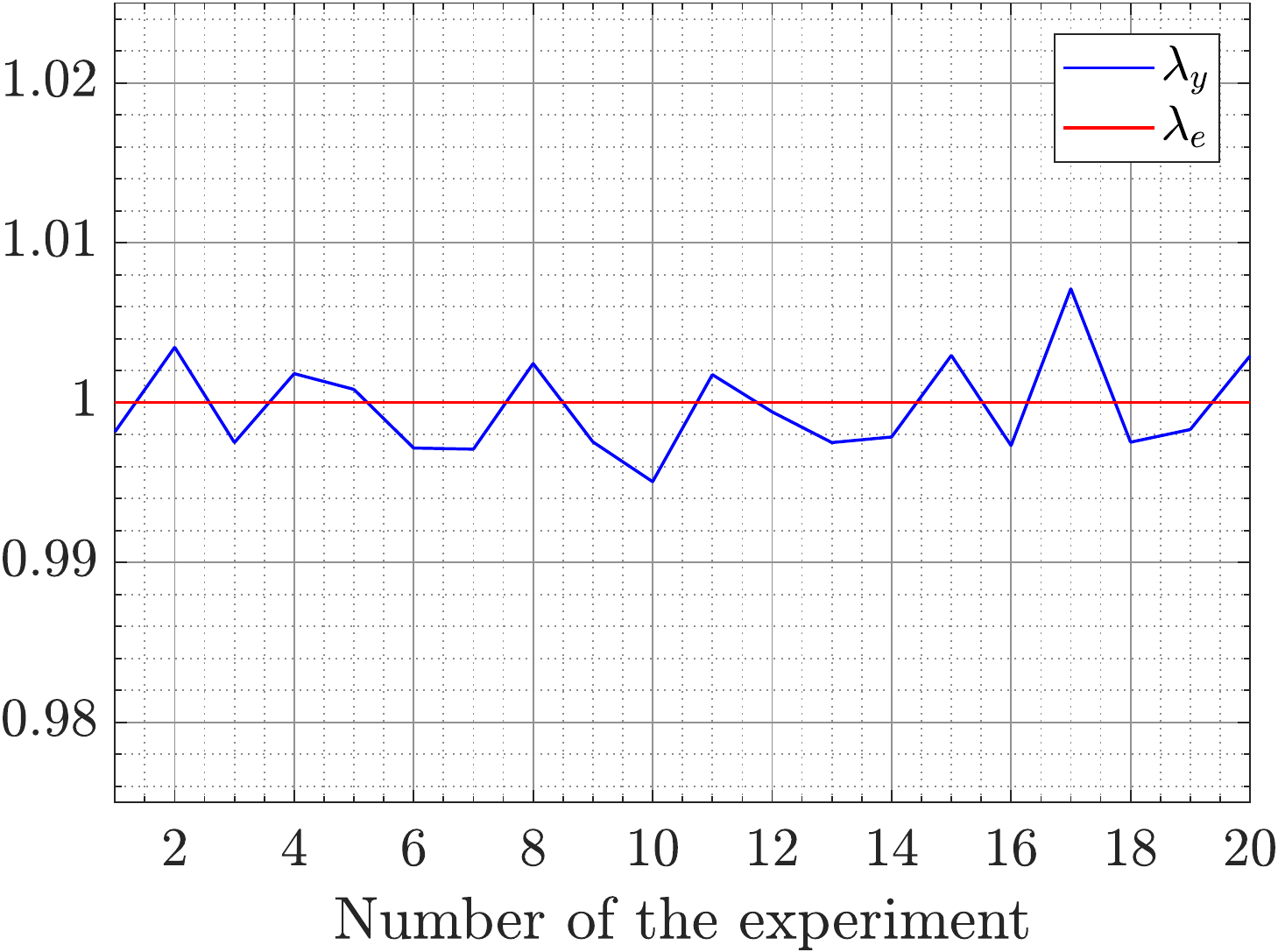}
}
\vspace{0.5cm}
\centerline{
\includegraphics[page=1, trim=0cm 0cm 0cm 0cm, clip, width=0.5\textwidth]{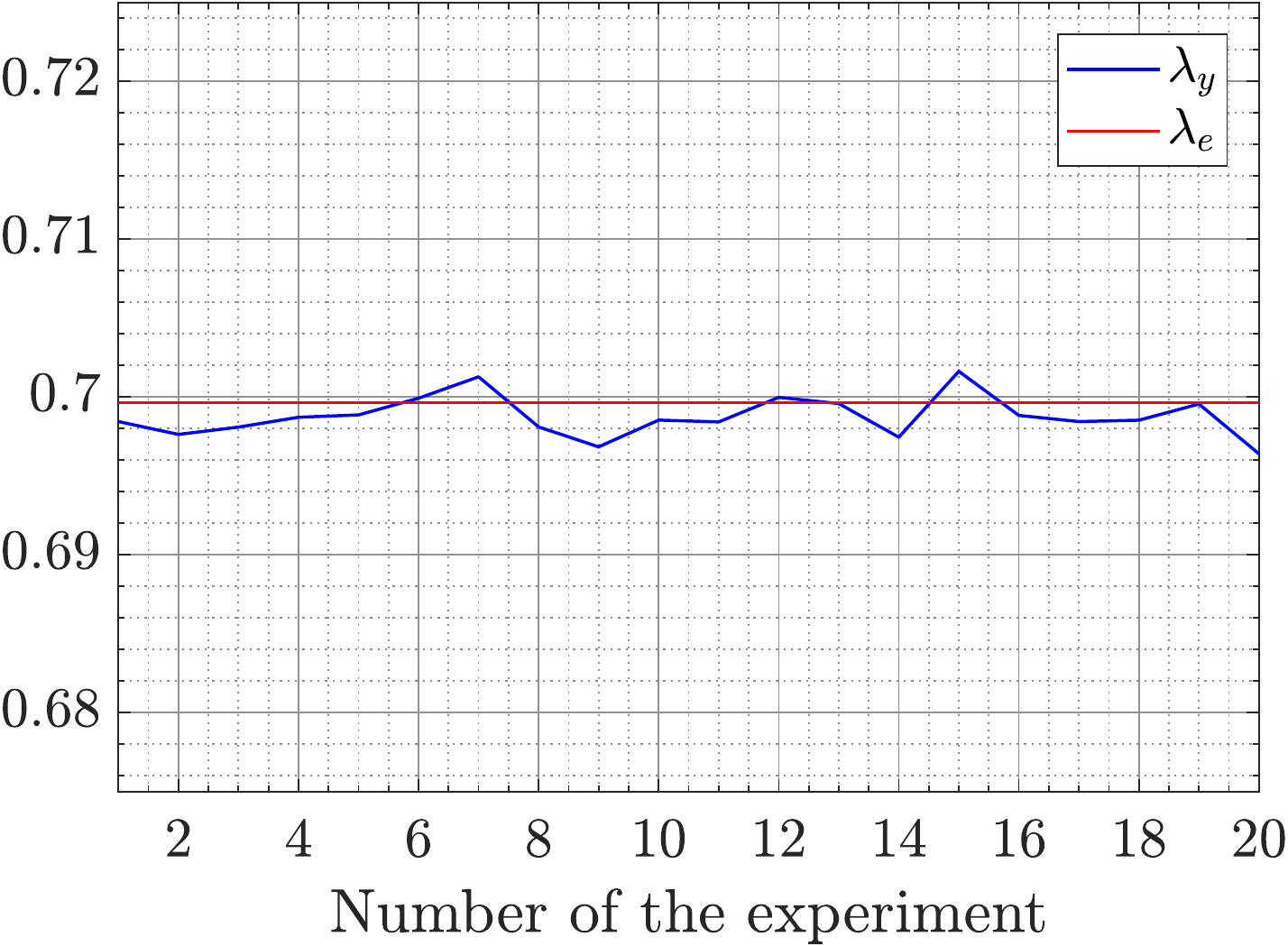}
\hfill
\includegraphics[page=1, trim=0cm 0cm 0cm 0cm, clip, width=0.5\textwidth]{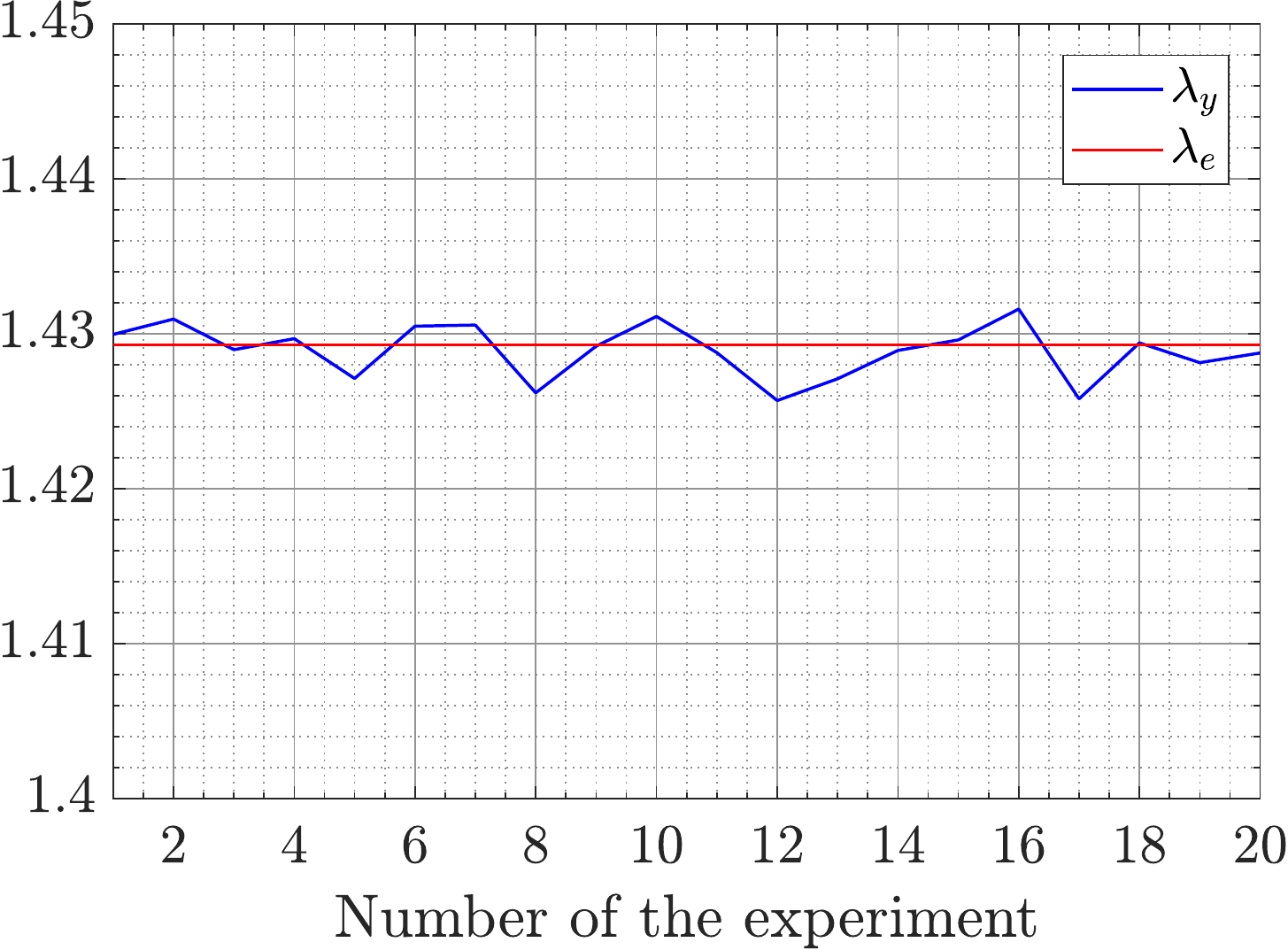}
}
\caption{The values of the effective thermal conductivity $\lambda_y$ and the estimated effective conductivity $\lambda_e$ in~\eqref{eq:CMA3} (for the domain $\Omega$ with $m=2000$ circular holes in Example~\ref{ex:2000NV}) vs. the number of the experiment for Case I (first row), Case II (second row, left), and Case III (second row, right).}
\label{fig:2000NVLam}
\end{figure}

\begin{example}\label{ex:2000NV}{\rm We take $m=2000$ circular holes of equal size with radius $r=0.0075$. In this example, the concentration $c=c(m,r) = 1000 r^2\pi\approx 0.1767$ is constant and the locations of these holes are chosen randomly. In this case, the following extension of CMA may be used 
	\begin{equation}\label{eq:CMA3}
		\lambda_e = \frac{1+c_1-c_2}{1-c_1+c_2} + O(c^3),
	\end{equation}
where $c_1$ denotes the conductor concentration, $c_2$ the insulator concentration, and $c=c_1+c_2$. 
Three cases are considered:
\begin{description}
	\item[Case I:] We assume that half of the holes are CNTs and the other half are voids (see Figure~\ref{fig:2000NV}). For this case, $c_1$ and $c_2$ are given by $c_1=c_2=500 r^2\pi\approx 0.0884$.
	\item[Case II:] All holes are voids, and hence $c_1=0$ while $c_2=1000 r^2\pi\approx 0.1767$.	
	\item[Case III:] All holes are CNTs, and hence $c_1=1000 r^2\pi\approx 0.1767$ while  $c_2=0$.
\end{description}

For each case, we run the code for 20 times, so that to get 20 different locations for these circular holes. In each of these 20 experiments, we compute the value of the effective thermal conductivity $\lambda_y$ by the presented method and the values of the estimated effective conductivity $\lambda_e$ by~\eqref{eq:CMA}.
As we can see from Figure~\ref{fig:2000NVLam}, $\lambda_e$ is a constant and the values of $\lambda_y$ depend on the locations of the holes. 
}
\end{example}

\begin{figure}[ht] %
\centerline{
\includegraphics[page=1, trim=0cm 0cm 0cm 0cm, clip, width=0.55\textwidth]{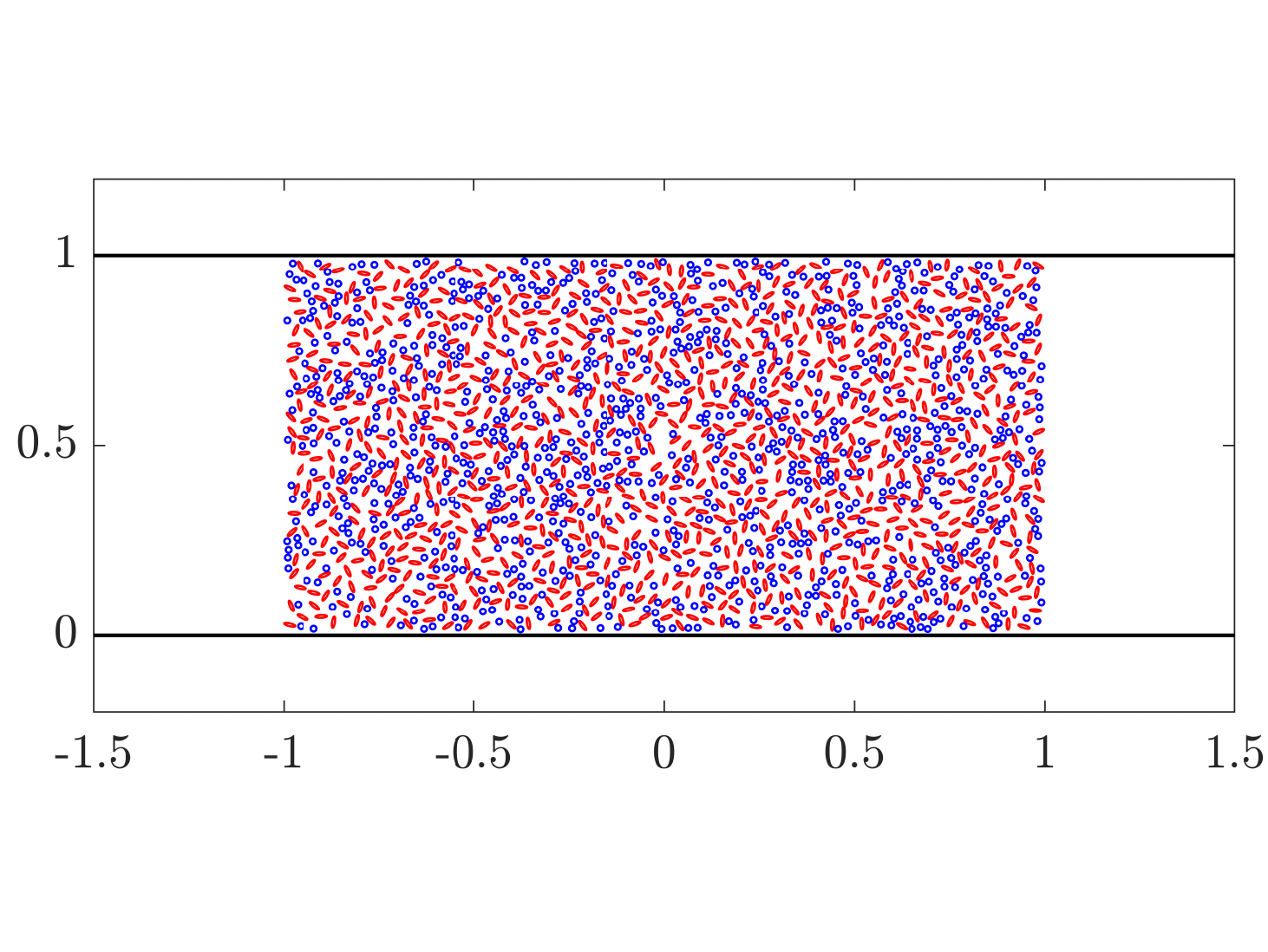}
\hfill
\includegraphics[page=1, trim=0cm 0cm 0cm 0cm, clip, width=0.5\textwidth]{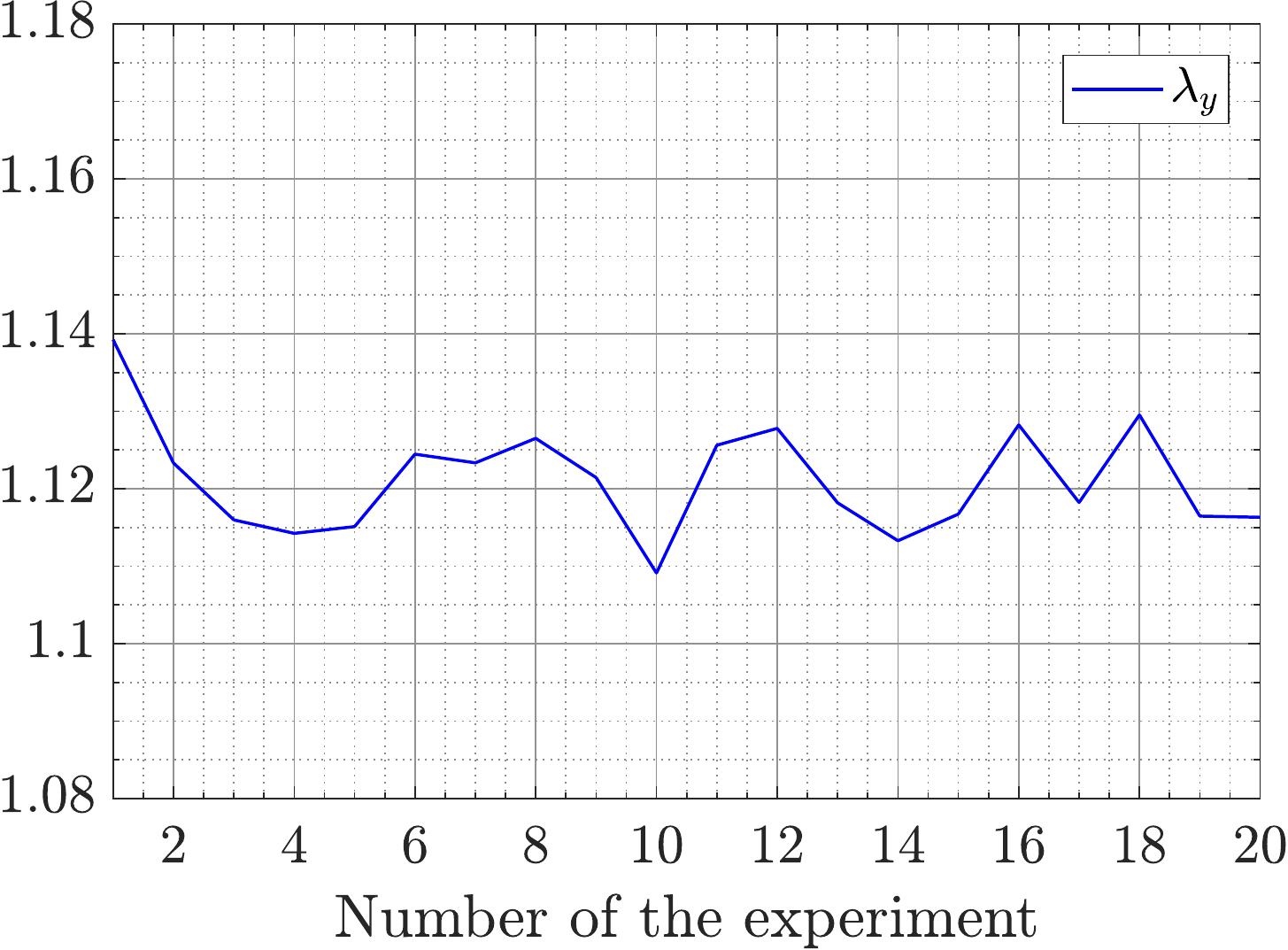}
}
\caption{On the left, the domain $\Omega$ in Example~\ref{ex:2000EC} with $m=p+\ell=2000$ holes, $p=1000$ circular voids (blue) and $\ell=1000$ elliptic CNTs (red). On the right, the values of the effective thermal conductivity $\lambda_y$ vs. the number of the experiment.}
\label{fig:2000EC}
\end{figure}

\begin{example}\label{ex:2000EC}{\rm 
We consider $m=2000$ elliptic and circular holes with $\ell=1000$ elliptic perfect conductors and $p=1000$ circular insulators of equal area $\pi r^2$ (see Figure~\ref{fig:2000EC}). The radius $r$ is chosen to be the same as in the previous example, i.e., $r=0.0075$. The locations of both elliptic and circular holes are chosen randomly. For the ellipses, we assume that the ratio between the length of the major axis and the minor axis is $4$, and the angles between the major axis and the $x$-axis are chosen randomly. As in the previous example, we run the code for 20 times. In each of these 20 experiments, we compute the value of the effective thermal conductivity $\lambda_y$ by the presented method. The computed values are shown in Figure~\ref{fig:2000EC} (left).

Since we have the same number of elliptic perfect conductors and circular insulators of equal area $\pi r^2$, the conductor concentration $c_1$ and the insulator concentration $c_2$ are equal and given by $c_1=c_2=500 r^2\pi\approx 0.0884$.
Although $c_1$ and $c_2$ here are the same as in Case I of the previous example, it is clear from Figures~\ref{fig:2000NVLam} (first row) and~\ref{fig:2000EC} (right) that the values $\lambda_y$ in this example (elliptic conductors) are larger than those in the previous example (circular conductors). 
}
\end{example}

\clearpage

\subsection{The dependence of $\lambda_y$ on $\phi$ and $c$}

We consider now a domain domain $\Omega$ with $m=\ell=276$ non-overlapping elliptic CNTs without any void. We assume that all CNTs are of equal size and elliptic shape. The ellipses are parametrized by~\eqref{eq:ellipse} with $a<b$, which means they are taken to be vertical. 

\begin{figure}[ht] %
\centerline{
\includegraphics[page=1, trim=0cm 2cm 0cm 2cm, clip, width=0.6\textwidth]{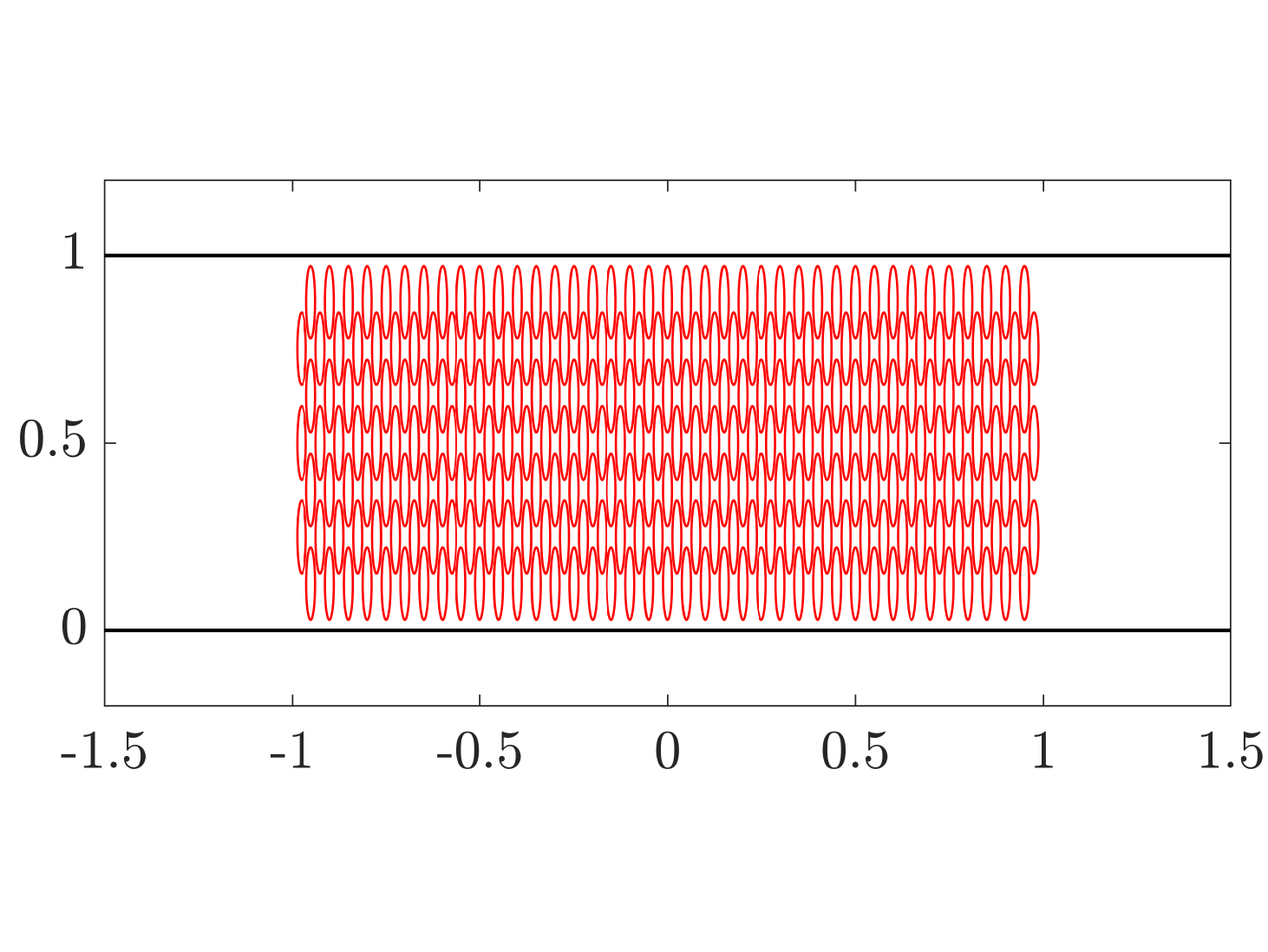}
}
\caption{The domain $\Omega$ with $m=276$ elliptic CNTs for $c=0.5$ and $\phi=1.3$.}
\label{fig:fig276}
\end{figure}

\begin{figure}[ht] %
\centerline{
\includegraphics[page=1, trim=0cm 0cm 0cm 0cm, clip, width=0.5\textwidth]{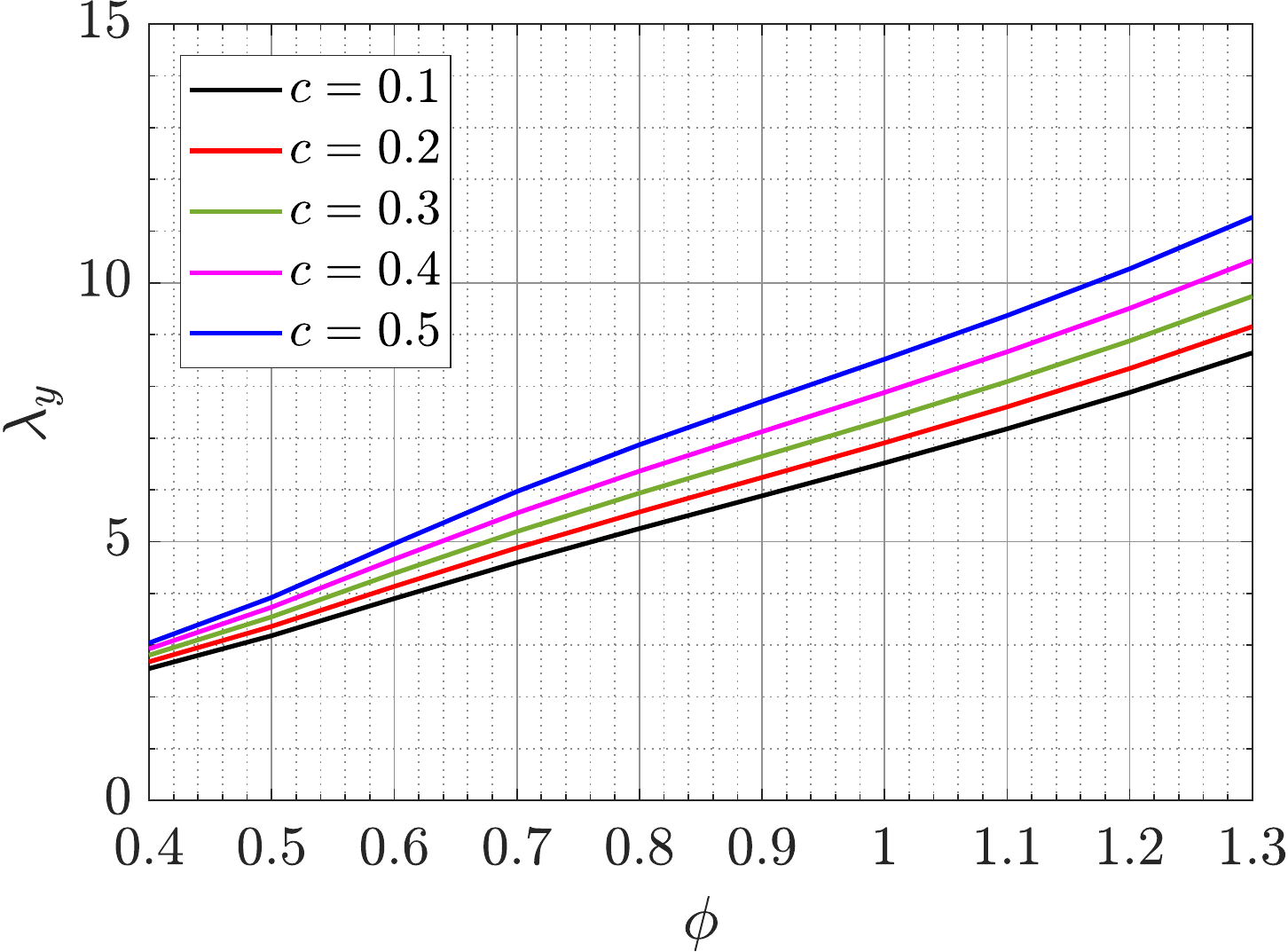}
\hfill
\includegraphics[page=1, trim=0cm 0cm 0cm 0cm, clip, width=0.5\textwidth]{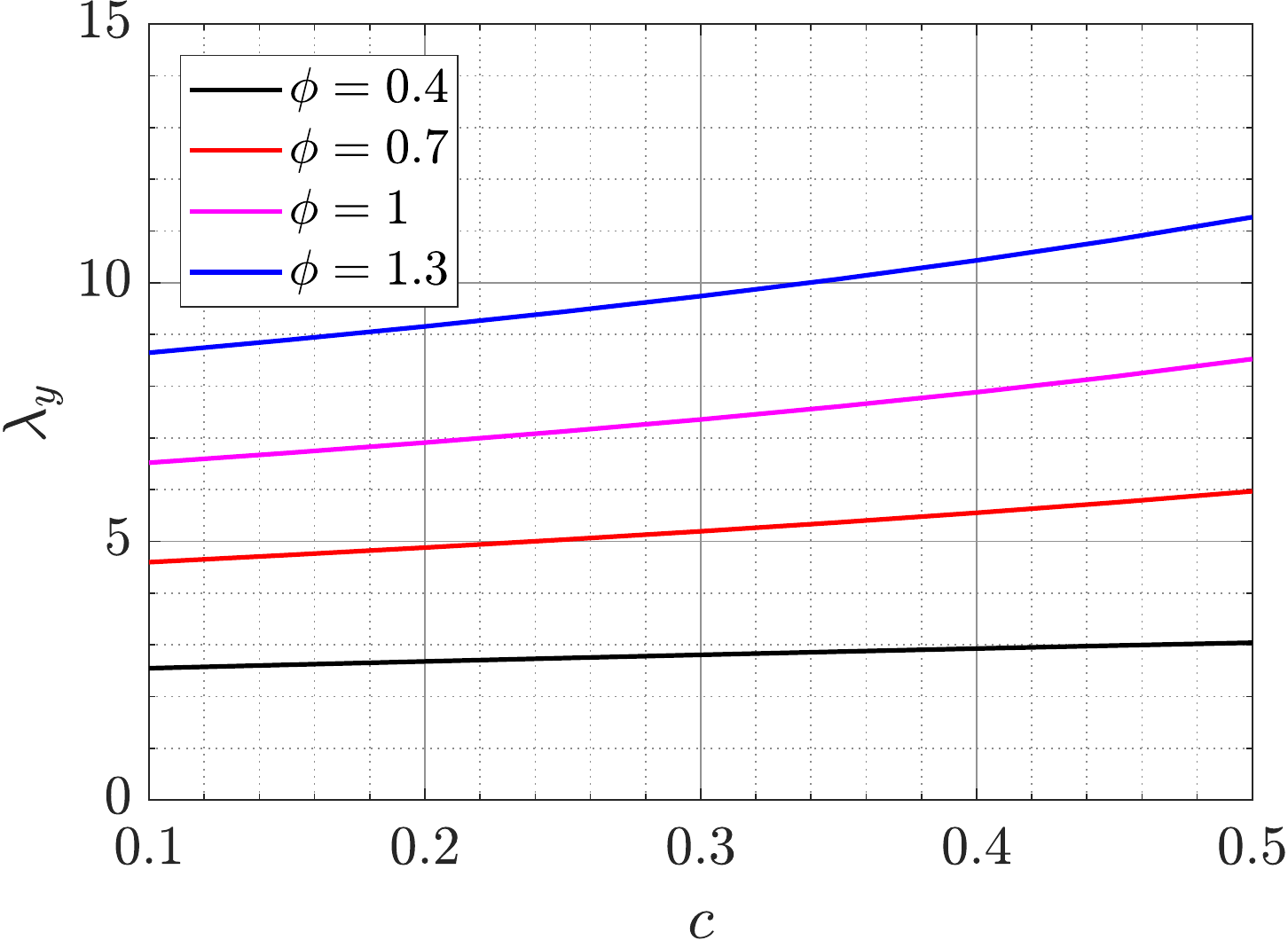}
}
\caption{The effective thermal conductivity $\lambda_y$ for the domain $\Omega$ with $m=276$ elliptic CNTs. On the left, the values of $\lambda_y=\lambda_y(\phi)$ for $\phi\in[0.4,1.3]$ and for several values of $c$. On the right, the values of $\lambda_y=\lambda_y(c)$ for $c\in[0.1,0.5]$ and for several values of $\phi$.}
\label{fig:Lam276}
\end{figure}

First we assume that the concentration $c=c(m,a,b)$ is constant and we choose the values of the parameters $a$ and $b$ such that the plane slits density $\phi=\phi(m,a,b)\in[0.4,1.3]$. The domain $\Omega$ for $c=0.5$ and $\phi=1.3$ is shown in Figure~\ref{fig:fig276}. 
We consider as well five values of  the concentration, $c=0.1,0.2,0.3,0.4,0.5$. Then, for each of these values, we compute and show in Figure~\ref{fig:Lam276} (left) the values of $\lambda_y=\lambda_y(\phi)$.

Afterwards, we take up the plane slits density $\phi=\phi(m,a,b)$ to be constant and we choose the values of the parameters $a$ and $b$ such that the concentration $c=c(m,a,b)\in[0.1,0.5]$. We consider four values of $\phi$, $\phi=0.4,0.7,1,1.3$, and compute again $\lambda_y=\lambda_y(c)$ for each case. The obtained results are presented in Figure~\ref{fig:Lam276} (right).

\section{Conclusion}
				
A systematic estimation of the local fields and the effective conductivity properties of 2D composites reinforced by uniformly and randomly distributed CNTs is carried out. It is assumed that the medium may contains voids as well. The CNTs are considered as perfectly conducting elliptic inclusions and the voids as circular insulators. For definiteness, a composite strip is considered with the given constant external field passing through the strip. The local field is governed by the Laplace equation in the multiply connected domain formed by the strip without two types of holes, CNTs and voids. The Dirichlet boundary condition is imposed on the CNTs boundary and the Neumann boundary condition governs the void boundary. A numerical method is developed to solve the mixed problem for a large number of CNTs and voids. The method is based on using the boundary integral equation with the generalized Neumann kernel~\cite{Weg-Nas,Nas-ETNA}. One key feature of this method is that it can be employed for domains with complex geometry as it provides accurate results even when the boundaries are close together. To solve the integral equation, the Fast Multipole Method has been employed, which enables to treat the case of thousands of CNTs and voids. With the help of conformal mappings, the presented method can be extended to include the case when CNTs and voids are rectilinear slits as done in~\cite{NasMo} for example.

The computational study has shown a dependence of the local fields and the effective conductivity $\lambda_y$ on the concentration of voids $c$ given by \eqref{eq:ell-conc} as well as on the density $\phi$ of CNTs given by \eqref{eq:ell-concP}. Besides the opposite conductive properties, voids and CNTs have also different types of the geometric parameters $c$ and $\phi$ that are not reduced to each other. Hence, the present study is concerned additionally with the case of three-phase composites of high contrast conductivity. It is demonstrated that our simulations are covered with the classical lower order approximations (Clausius-Mossotti, Maxwell) for dilute composites. The high order concentrations and densities led to different results from the classical ones. It is worth noting that the huge number of numerical experiments for uniformly distributed inclusions yield the graphical dependencies of $\lambda_y$ on $c$ and $\phi$, which can be used in practical applications.

\section*{Acknowledgments}

\bibliographystyle{unsrt}
\bibliography{references}

\end{document}